\documentclass[12pt]{article}

\usepackage{amsthm}
\usepackage[dvips]{graphicx}

\usepackage{epsfig}

\usepackage{latexsym}
\usepackage{amsmath}
\usepackage{amssymb}
\usepackage{amsfonts}

\def\no{\noindent}
\theoremstyle{definition}

\newtheorem{dfn}{Definition}[section]

\newtheorem{ex}[dfn]{Example}
\newtheorem{rem}[dfn]{Remark}

\theoremstyle{plain}

\newtheorem{add}[dfn]{Addendum}
\newtheorem{thm}[dfn]{Theorem}

\newtheorem{lem}[dfn]{Lemma}

\newtheorem{prop}[dfn]{Proposition}

\newtheorem{cor}[dfn]{Corollary}

\def\proof{\par\medskip\noindent{\it Proof. }}

\def\lra{\longrightarrow}

\def\ra{\rightarrow}
\def\Ra{\Rightarrow}

\def\half{\frac{1}{2}}
\def\C{{\Bbb C}}
\def\R{{\Bbb R}}

\def\Z{{\Bbb Z}}

\def\P{{\Bbb P}}

\def\N{{\Bbb N}}

\def\eps{\epsilon}
\def\al{\alpha}

\def\ga{\gamma}

\def\de{\delta}
\def\De{\Delta}
\def\Si{\Sigma}
\def\si{\sigma}

\def\la{\lambda}
\def\La{\Lambda}
\def\Om{\Omega}
\def\om{\omega}

\def\acts{\curvearrowright}
\def\D{\partial}

\def\embed{\hookrightarrow}
\def\8{\infty}
\def\<{\langle}
\def\>{\rangle}
\def\geo{\partial_{\infty}}
\def\tits{\partial_{Tits}}
\def\tangle{\angle_{Tits}}

\def\ol{\overline}

\newcommand{\restr}{\mbox{\Large \(|\)\normalsize}}
\def\3{\ss}

\def\goth{\mathfrak}

\begin{document}

\title{Convex functions on symmetric spaces,
side lengths of polygons 
and the stability inequalities for weighted configurations at infinity} 
\author{Misha Kapovich, Bernhard Leeb and John Millson}
\date{November 6, 2005}

\maketitle

\begin{abstract}
\no
In a symmetric space of noncompact type $X=G/K$ oriented geodesic segments 
correspond modulo isometries 
to vectors in the Euclidean Weyl chamber. 
We can hence assign vector valued lengths to segments. 
Our main result is a system of
homogeneous linear inequalities, 
which we call the generalized triangle inequalities 
or stability inequalities, 
describing the restrictions on the 
vector valued side lengths of oriented polygons. 
It is based on the mod $2$ Schubert calculus in the real Grassmannians $G/P$
for maximal parabolic subgroups $P$. 

The side lengths of polygons in Euclidean buildings are studied in 
the related paper \cite{polbuil}.
Applications of the geometric results in both papers 
to algebraic group theory are given in \cite{alg}.
\end{abstract}

\tableofcontents

\section{Introduction} 

When studying asymptotic properties of the spectra 
of certain linear partial differential equations 
in mathematical physics 
Hermann Weyl was led in \cite{Weyl}
to the question 
how the spectra of two compact self adjoint operators 
are related to the spectrum of their sum. 
The restrictions turn out to be homogeneous linear inequalities 
involving finite subsets of eigenvalues.
It suffices to understand the question in the finite-dimensional case
where it can be phrased as follows. 
Here 
$\al = (\al_1,\dots,\al_m)$, 
$\beta = (\beta_1,\dots,\beta_m)$ and
$\ga = (\ga_1,\dots,\ga_m)$ 
denote $m$-tuples of real 
numbers arranged in decreasing order
and with sum equal to zero.

\medskip\no
{\bf Eigenvalues of a sum problem.}
{\em Give necessary and sufficient conditions on $\al$, $\beta$ and $\ga$ 
in order that there exist traceless Hermitian matrices 
$A,B,C\in i\cdot su(m)$
with spectra $\al,\beta,\ga$
and satisfying}
\begin{equation*}
A+B+C=0.
\end{equation*}
There is a multiplicative version of this question. 
We recall that the {\em singular values} of a matrix $A$ in $GL(m,\C)$ 
are defined as the (positive) square roots 
of the eigenvalues of the matrix $AA^*$. 

\medskip\no
{\bf Singular values of a product problem.}
{\em Give necessary and sufficient conditions on $\al$, $\beta$ and $\ga$ 
in order that there exist matrices 
$A,B,C\in SL(m,\C)$
the logarithms of whose singular values are $\al,\beta,\ga$
and which satisfy}
\begin{equation*} 
ABC=1 .
\end{equation*}

\no
We refer to \cite{Fulton2000} 
for detailed information on these questions and their history. 

Both questions have natural geometric interpretations and generalizations.
Let us consider the group $G=SL(m,\C)$,
its maximal compact subgroup $K=SU(m)$ 
and the symmetric space $X=G/K$. 
Decompose the Lie algebra  ${\goth g}=sl(m,\C)$ of $G$ according to
${\goth g} = {\goth k} \oplus {\goth p}$ where ${\goth k}=su(m)$ 
is the Lie algebra of $K$
and ${\goth p}=i\cdot su(m)$ is the orthogonal complement of ${\goth k}$ relative to the 
Killing form. 
${\goth p}$ is canonically identified with the tangent space $T_oX$ to $X$ 
at the base point $o$ fixed by $K$.

The singular values of a group element $g\in G$
form a vector $\si(g)$ in the Euclidean Weyl chamber $\De_{euc}$
and are a complete invariant of the double coset $KgK$. 
More geometrically, they can be interpreted as a vector valued distance:
Given two points $x_1,x_2\in X$, $x_i=g_io$, 
the singular values $\si(g_1^{-1}g_2)$ 
are the complete invariant of the pair $(x_1,x_2)$
modulo the $G$-action.
We will call it the {\em $\De$-length} of the oriented geodesic segment $\ol{x_1x_2}$. 
The Singular values of a product problem
thus asks about the possible $\De$-side lengths of triangles in $X$.

In the same vein,
the Eigenvalues of a sum problem is a problem for triangles in ${\goth p}$
equipped with the geometry, in the sense of Felix Klein,
having as automorphisms 
the group $Aff({\goth p})$ 
generated by the adjoint action of $K$ and all translations. 
We call $({\goth p},Aff({\goth p}))$ 
the {\em infinitesimal symmetric space} associated to $X$. 
In this geometry, 
pairs of points are equivalent 
if and only if their difference matrices have equal spectra.
The spectra of the matrices can again be interpreted as $\De$-lengths,
and Weyl's question amounts to finding the 
sharp {\em triangle inequalities} 
in this geometry.

This paper is devoted to the 

\medskip\no
{\bf Problem.}
{\em Study for an arbitrary connected semisimple real Lie group $G$
of noncompact type with finite center 
the spaces ${\cal P}_n(X),{\cal P}_n({\goth p})\subset\De_{euc}^n$,
$n\geq3$, 
of possible $\De$-side lengths of oriented $n$-gons 
in the associated Riemannian symmetric space $X$
and the corresponding infinitesimal symmetric space ${\goth p}$.}

\medskip
Our main result is an explicit description 
of ${\cal P}_n(X)$ and ${\cal P}_n({\goth p})$ 
in terms of a finite system of homogeneous linear inequalities 
parametrized by the Schubert calculus associated to $G$,
see Theorem \ref{generalschubertinequalitiesintro} below. 
In particular,
both spaces are finite-sided polyhedral cones 
and we will refer to them 
as {\em $\De$-side length polyhedra}. 
Our approach is based on differential-geometric techniques 
from the theory of nonpositively curved spaces.

To determine the side length polyhedra
we relate oriented polygons in the symmetric space $X$, respectively,
the infinitesimal symmetric space ${\goth p}$ 
via a Gauss map type construction 
to weighted configurations 
on the spherical Tits building at infinity $\tits X$. 
The question on the possible $\De$-side lengths of polygons 
then translates into a question about weighted configurations 
as it occurs in geometric invariant theory.
Namely, 
after suitably generalizing the concept of Mumford stability, 
it turns out that 
the set of possible $\De$-side lengths of polygons 
coincides with the set of possible $\De$-weights 
(defined below) 
of semistable configurations 
(Theorems \ref{equivalenceofmoduli} and \ref{polxconf:samelengths}),
that is,
we have to determine the $\De$-weights 
for which there exist semistable configurations.

Since the questions for polygons in $X$ respectively ${\goth p}$
translate into the same question for configurations,
we obtain as a byproduct: 

\begin{thm}
\label{thm:thompson}
${\cal P}_n(X)={\cal P}_n({\goth p})$.
\end{thm}
This 
generalizes the Thompson Conjecture \cite{Thompson}
which was proven for $GL(m,\C)$ in \cite{Klyachko2}
and more generally for complex semisimple groups in 
\cite{AMW}. 
Another proof of the most general
version of the Thompson conjecture has recently been given in 
\cite{EvensLu}.

In our (logically independent) paper \cite{polbuil} 
we investigate 
the $\De$-side lengths of polygons in {\em Euclidean buildings}.
The main result asserts that 
for a thick Euclidean building $Y$ 
the $\De$-side length space ${\cal P}_n(Y)$ 
depends only on the spherical Coxeter complex associated to $Y$.
The proof exploits an analogous relation between polygons in $Y$ 
and weighted configurations on the Tits boundary $\tits Y$ 
by ways of a Gauss map.
Along the way we show that 
${\cal P}_n(Y)$ coincides with the space of $\De$-weights 
of semistable weighted configurations on the spherical building $\tits Y$
and, moreover, 
that the space of $\De$-weights 
of semistable weighted configurations on a thick spherical building $B$
depends only on the spherical Coxeter complex attached to $B$.
Note that every spherical Tits building occurs as the Tits boundary 
of a Euclidean building, 
for instance, of the complete Euclidean cone over itself. 

The results of \cite{polbuil} imply 
for polygons in Riemannian symmetric spaces: 

\begin{thm} 
\label{polddepsphcox}
${\cal P}_n(X)$ 
depends only on the spherical Coxeter complex attached to $X$. 
\end{thm}

In \cite{alg} we apply the results of this paper and \cite{polbuil}
to algebra. 
The generalized triangle inequalities give necessary conditions for solving
a number of problems in algebraic group theory. 
For the case $G=GL(m)$ 
we give a new proof of the Saturation Conjecture 
first proved in \cite{KnutsonTao}. 

\medskip
Most of the remaining part of the introduction 
will be devoted to describing 
the (semi)stability inequalities 
for $\De$-weights of configurations on $\geo X$.
As we said earlier, 
they coincide with the inequalities 
for the $\De$-side lengths of polygons in $X$. 

A symmetric space of noncompact type $X$  
is a complete simply connected Riemannian manifold 
with nonpositive sectional curvature
and as such can be compactified to a closed ball by attaching 
an ideal boundary (sphere) $\geo X$.
This construction generalizes the compactification of hyperbolic space
given by the conformal Poincar\'e ball model.
The natural $G$-action on $X$ by isometries 
extends to a continuous action on $\geo X$. 
The $G$-orbits on $\geo X$ are parametrized by the spherical Weyl chamber,
$\geo X/G\cong\De_{sph}$. 
They are homogeneous $G$-spaces of the form $G/P$ with $P$ a parabolic
subgroup
and we call them {\em generalized flag manifolds}. 

A {\em weighted configuration} on $\geo X$ 
is a map 
$\psi:(\Z/n\Z,\nu)\ra \geo X$ 
from a finite measure space.
We think of the masses $m_i:=\nu(i)$ 
placed in the points $\xi_i:=\psi(i)$ at infinity. 
The {\em $\De$-weights} 
$h=(h_1,\dots,h_n)\in\De_{euc}^n$ 
of the configuration 
contain the information on the masses and the $G$-orbits where they are
located:
To each orbit $G\xi_i$ corresponds the point $acc(\xi_i)$ in $\De_{sph}$
where $acc:\geo X\ra\De_{sph}$ denotes the natural projection. 
We view the spherical simplex $\De_{sph}$ 
as the set of unit vectors in the complete Euclidean cone $\De_{euc}$ 
and define 
$h_i:=m_i\cdot acc(\xi_i)$.

To a weighted configuration on $\geo X$
one can associate a natural {\em convex} function on $X$,
the {\em weighted Busemann function}
$\sum_i m_i\cdot b_{\xi_i}$ 
(well-defined up to an additive constant),
compare \cite{DouadyEarle}. 
The Busemann function $b_{\xi_i}$
measures the relative distance from the point $\xi_i$ at infinity.
We define {\em stability} and {\em semistability}
of a weighted configuration 
in terms of asymptotic properties of its Busemann function. 
These asymptotic properties can in fact 
be expressed in terms of the Tits angle metric on $\geo X$ 
which leads to a notion of stability for weighted configurations
on abstract spherical buildings,
see also \cite{polbuil}. 
Our notion of stability agrees with Mumford stability 
from geometric invariant theory 
when $G$ is a complex group. 
Examples can be found in section \ref{sec:expl}
where we determine explicitly 
the (semi)stable weighted configurations 
(more generally, of finite measures)
on the Grassmannians associated to the classical groups. 

The possible $\De$-weights 
for semistable weighted configurations on $\geo X$ 
are given by 
a finite system of homogeneous linear inequalities. 
We first describe the {\em structure} of these inequalities.
Let $(S,W)$ denote the spherical Coxeter complex attached to $G$ 
and think of the spherical Weyl chamber $\De_{sph}$ as being embedded in $S$. 
For every vertex $\zeta$ of $\De_{sph}$ 
and every $n$-tuple of vertices $\eta_1,\dots,\eta_n\in W\zeta$ 
we consider the inequality
\begin{equation}
\label{stabineq:intro}
\sum_i m_i\cdot\cos\angle(\tau_i,\eta_i)\leq0 ,
\end{equation}
for $m_i\in\R^+_0$ and $\tau_i\in\De_{sph}$
where $\angle$ measures the spherical distance in $S$.
We may rewrite the inequality as follows 
using standard terminology of Lie theory: 
Let $\la_{\zeta}\in\De_{euc}$ be the fundamental coweight 
contained in the edge with direction $\zeta$,
and let $\la_i:=w_i\la_{\zeta}$ where $[w_i]\in W/W_{\zeta}$ 
such that $w_i\zeta=\eta_i$. 
With the renaming $h_i=m_i\tau_i$ of the variables 
(\ref{stabineq:intro}) becomes the homogeneous linear inequality
\begin{equation}
\label{stabineqlie:intro}
\sum_i \langle h_i,\la_i\rangle \leq0 .
\end{equation}
The full family of these inequalities
has only the trivial solution. 
The stability inequalities are given by a subset
of these inequalities
which we single out using the {\em Schubert calculus}.

For a vertex $\zeta$ of $\De_{sph}$ 
we denote by $Grass_{\zeta}\subset\geo X$ 
the corresponding maximally singular $G$-orbit on $\geo X$. 
We call it a {\em generalized Grassmannian}
because in the case of $SL(m)$ 
the $Grass_{\zeta}$ are the usual Grassmann manifolds. 
The stabilizers of points in $Grass_{\zeta}$ 
are the conjugates of a maximal parabolic subgroup 
$P\subset G$.
The restriction of the $G$-action to $P$
stratifies $Grass_{\zeta}$ into {\em Schubert cells},
one cell $C_{\eta_i}$ corresponding to each vertex $\eta_i\in S$ 
in the orbit $W\zeta$ of $\zeta$ under the Weyl group $W$. 
Hence, if we denote $W_{\zeta}:=Stab_W(\zeta)$ 
then the Schubert cells correspond to cosets in $W/W_{\zeta}$. 
The {\em Schubert cycles} are defined as the closures $\ol C_{\eta_i}$ 
of the Schubert cells;
they are unions of Schubert cells.
As real algebraic varieties, 
the Schubert cycles represent homology classes 
$[\ol C_{\eta_i}]\in H^{\ast}(Grass_{\zeta};\Z/2\Z)$
which we abbreviate to $[C_{\eta_i}]$;
in the complex case they even represent {\em integral} homology classes. 

Now we can formulate our main result. 
We recall that ${\cal P}_n(X)={\cal P}_n({\goth p})$
coincides with the set of $\De$-weights 
of semistable weighted configurations on $\geo X$. 

\begin{thm}
[Stability inequalities for noncompact semisimple Lie groups]
\label{generalschubertinequalitiesintro}

(i)
The set ${\cal P}_n(X)$ consists of all $h\in\De_{euc}^n$ 
such that 
(\ref{stabineqlie:intro})
holds whenever the intersection 
of the Schubert classes 
$[C_{\eta_1}],\dots,[C_{\eta_n}]$ 
in $H_{\ast}(Grass_{\zeta};\Z/2\Z)$ 
equals $[pt]$.

(ii)
If $G$ is complex,
then 
the set ${\cal P}_n(X)$ consists of all $h\in\De_{euc}^n$ 
such that 
the inequality 
(\ref{stabineqlie:intro})
holds whenever the intersection 
of the integral Schubert classes 
$[C_{\eta_1}],\dots,[C_{\eta_n}]$ 
in $H_{\ast}(Grass_{\zeta};\Z)$ 
equals $[pt]$.
\end{thm}

Our argument shows moreover that 
the system obtained by imposing all inequalities 
(\ref{stabineqlie:intro}) 
whenever the intersection 
of the Schubert classes 
$[C_{\eta_1}],\dots,[C_{\eta_n}]$ is {\em nonzero} in  
$H_{\ast}(Grass_{\zeta};\Z/2\Z)$ has the same set of solutions as the
smaller system obtained when we require the intersection to be $[pt]$. 
Note that
in the complex case 
the set of necessary and sufficient inequalities 
obtained from the {\em integral Schubert} calculus is in general smaller.

Interestingly,
the stability inequalities 
given by Theorem \ref{generalschubertinequalitiesintro} 
depend on the Schubert calculus
whereas their solution set
depends only on the Weyl group.
This is due to possible redundancies. 

In rank one,
i.e.\ when $X$ has strictly negative sectional curvature,
we have $\De_{euc}\cong\R^+_0$
and the stability inequalities are just the ordinary triangle inequalities. 
In section \ref{sec:rank2} 
we determine the side length spaces for all 
symmetric spaces of rank two. 
The rank three case is already quite involved 
and it is treated in the paper \cite{KumarLeebMillson}. 

The polyhedron ${\cal P}_3({\goth p})$
was first determined for $G=SL(m,\C)$ 
by Klyachko \cite{Klyachko} who proved that the
inequalities corresponding to triples of Schubert classes with
intersection a positive multiple of the point class were
{\em sufficient}. The necessity of these inequalities has been
know for some time, see \cite[sec.\ 6, Prop.\ 2]{Fulton2000} for a proof of
their necessity and the history of this proof.
Klyachko's theorem was refined by Belkale \cite{Belkale} 
who showed that it suffices to restrict
to those triples of Schubert classes whose intersection is the
point class. The determination of ${\cal P}_3({\goth p})$ for
general complex simple $G$
was accomplished in
\cite{BerensteinSjamaar}
using methods from algebraic geometry.
However they gave a larger system (than ours) 
consisting of all inequalities
where the intersections of Schubert
classes is a {\em nonzero multiple} of the point class. 
Thus our Theorem \ref{generalschubertinequalitiesintro} 
for the complex case 
is a refinement of their result.
In the general real case the polyhedron ${\cal P}_3({\goth p})$
was determined in \cite{OsheaSjamaar}.
However their inequalities are quite different from ours.
They are associated to the
{\em integral} Schubert calculus of the {\em complexification}
${\goth g}\otimes {\goth C}$
and are efficient for the case of split ${\goth g}$
but become less and less efficient
as the real rank of ${\goth g}$ (i.e.\ the rank of the symmetric space $X$)
decreases.
For instance, for the case of real rank one
they have a very large number of inequalities when the
ordinary triangle inequalities alone will suffice.
This is recognized in \cite{OsheaSjamaar}
and the problem is posed
(Problem 9.5 on page 451)
as to whether a formula of the
type we found above in terms of the Schubert calculus {\em modulo $2$}
would exist.

\medskip
The paper is organized as follows. 
In section \ref{sec:geomprelim}
we provide some background 
from the geometry of spaces of nonpositive curvature, 
symmetric spaces of noncompact type
and of spherical buildings. 
In section \ref{sec:conf} 
we define and study a notion of stability 
for measures and weighted configurations 
on the ideal boundary of symmetric spaces of noncompact type. 
We provide analogues of some basic results in geometric invariant theory,
such as a Harder-Narasimhan Lemma 
(Theorem \ref{hardernarasimhan})
and prove our main result Theorem \ref{generalschubertinequalitiesintro}. 
The {\em weak} stability inequalities 
considered in section \ref{sec:wtrineq}
correspond to particularly simple intersections of Schubert cycles
and they have a beautiful geometric interpretation in terms of convex hulls. 
In section \ref{git}
we explain in the example 
of weighted configurations on complex projective space 
that our notion of stability 
matches with Mumford stability.
In section \ref{sec:pol}
we discuss the relation between 
polygons in $X$ and configurations on $\geo X$
and prove the generalized Thompson Conjecture (Theorem \ref{thm:thompson}). 
In section \ref{sec:expl}
we will make explicit 
the stability condition for measures 
supported on the (generalized) Grassmannians associated to the classical
groups. 
In section \ref{sec:rank2} 
we make a detailed study of the polyhedra
${\cal P}_n({\goth p})$ 
for the rank $2$ complex simple groups.
We make our system of stability inequalities explicit 
and for $n=3$ we describe the minimal subsystems, 
i.e.\ the facets of the polyhedron.
Moreover we give the generators (edges) of ${\cal P}_3({\goth p})$.
The inequalities for the group $G_2$ were computed previously 
in \cite[Example 5.2.2]{BerensteinSjamaar}. 
The paper \cite{KumarLeebMillson} studies 
${\cal P}_3({\goth p})$ 
for the rank $3$ examples and describes the minimal subsystems
and the generators of the cone for the root systems $C_3$ and $B_3$.

\medskip
{\em Acknowledgements.}
We would like to thank 
Andreas Balser for reading an early version of this paper
and Chris Woodward for helpful suggestions
concerning the computations in section \ref{sec:rank2}. 
Also we took the multiplication
table for the Schubert classes for $G_2$ from \cite{TelemanWoodward}. 
In section \ref{sec:rank2} 
we used the
computer program Porta written by Thomas Christof and Andreas L\"obel to
find the minimal subsystems and the generators of the cones. Finally
we would like to thank George Stantchev for finding the computer program 
Porta and for much help and advice in implementing it.

\section{Preliminaries}
\label{sec:geomprelim}

In this section,
mostly to fix our notation, 
we will briefly review some basic facts
about spaces of nonpositive curvature 
and especially Riemannian symmetric spaces of noncompact type.
We will omit most of the proofs.
For more details on spaces with upper curvature bound 
and in particular spaces with nonpositive curvature,
we refer to \cite[ch.\ 1-2]{Ballmann}, \cite[ch.\ 4+9]{Burago}, 
\cite[ch.\ 2]{KleinerLeeb} and \cite[ch.\ 2]{habil},
for the geometry of symmetric spaces of noncompact type 
to 
\cite{karpelevich}, 
\cite[ch.\ II.10]{BridsonHaefliger}, 
\cite[ch.\ 2+3]{Eberlein} and \cite[ch.\ 6]{Helgason}, 
and for the theory of spherical buildings from a geometric viewpoint,
i.e.\ within the framework of spaces with curvature bounded above,
to \cite[ch.\ 3]{KleinerLeeb}.

\subsection{Metric spaces with upper curvature bounds}

Consider a complete geodesic space,
that is,
a complete metric space $Y$
such that any two points $y_1,y_2\in Y$ 
can be joined by a rectifiable curve with length $d(y_1,y_2)$;
such curves are called {\em geodesic segments}. 
Note that we do not assume $Y$ to be locally compact.
Although there is no smooth structure nor a Riemann curvature tensor
around,
one can still make sense of a sectional curvature {\em bound}
in terms of distance comparison. 
We will only be interested in {\em upper} curvature bounds.
One says that $Y$ has (globally) curvature $\leq k$ 
if all triangles in $Y$ are thinner than corresponding triangles 
in the model plane (or sphere) $M^2_k$ of constant curvature~$k$.
Here,
a geodesic triangle $\De$ in $Y$ is a one-dimensional object 
consisting of three points and geodesic segments joining them. 
A comparison triangle $\tilde\De$ for $\De$ in $M^2_k$
is a triangle with the same side lengths.
To every point $p$ on $\De$ corresponds a point $\tilde p$ on
$\tilde\De$,
and we say that $\De$ is thinner than $\tilde\De$ 
if for any points $p$ and $q$ on $\De$ we have 
$d(p,q)\leq d(\tilde p,\tilde q)$. 
A metric space with curvature $\leq k$ is also called a {\em $CAT(k)$-space}.
It is a direct consequence of the definition 
that any two points are connected by a {\em unique} geodesic 
if $k\leq0$,
or if $k>0$ and the points have distance $<\frac{\pi}{\sqrt{k}}$.

By Toponogov's Theorem 
a complete simply-connected manifold has curvature $\leq k$ in the
distance comparison sense
if and only if it has sectional curvature $\leq k$.

\medskip
The presence of a curvature bound allows 
to define {\em angles} between segments
$\si_i:[0,\eps)\ra Y$ initiating in the same point $y=\si_1(0)=\si_2(0)$
and parametrized by unit speed.
Let $\tilde\al(t)$ be the angle 
of a comparison triangle for $\De(y,\si_1(t),\si_2(t))$
in the appropriate model plane 
at the vertex corresponding to $y$.
If $Y$ has an upper curvature bound 
then the comparison angle $\tilde\al(t)$ 
is monotonically decreasing 
as $t\searrow0$. 
It therefore converges,
and we define the angle $\angle_y(\si_1,\si_2)$ of the segments at $y$ 
as the limit.
In this way, 
one obtains a pseudo-metric 
on the space of segments emanating from a point $y\in Y$.
Identification of segments with angle zero and metric completion 
yields the {\em space of directions} $\Si_yY$.
One can show that if $Y$ has an upper curvature bound,
then $\Si_yY$ has curvature $\leq1$.

\subsection{Spaces of nonpositive curvature}
\label{sec:hadamard}

In this section, we assume 
that $Y$ is a space of nonpositive curvature.
We will call spaces of nonpositive curvature
also {\em Hadamard spaces}. 

A basic consequence of nonpositve curvature is
that the distance function $d:Y\times Y\ra\R^+_0$ is {\em convex}.
It follows that geodesic segments between any two points are unique
and globally minimizing.
In particular, $Y$ is contractible. 
We will call a complete geodesic $l\subset Y$ also a {\em line}
since it is an isometrically embedded copy of $\R$. 

\medskip
A (parametrized) geodesic {\em ray} is an isometric embedding
$\rho:[0,\infty)\lra Y$. 
Two rays $\rho_1$ and $\rho_2$ are called {\em asymptotic}
if $t\mapsto d(\rho_1(t),\rho_2(t))$ stays bounded 
and hence, by convexity, (weakly) decreases. 
An equivalence class of asymptotic rays is called an {\em ideal} point
or a point {\em at infinity}.
If a ray $\rho$ represents an ideal point $\xi$,
we also say that $\rho$ is {\em asymptotic} to $\xi$.
We define the {\em geometric boundary} $\geo Y$
as the set of ideal points.
The topology on $Y$ can be canonically extended to the 
{\em cone topology} on $\bar Y:=Y\cup\geo Y$
If $Y$ is locally compact, $\bar Y$ is a compactification of $Y$.

There is a natural metric on $\geo Y$,
the {\em Tits metric}. 
The Tits distance of two ideal points $\xi_1,\xi_2\in\geo Y$ is defined as
$\tangle(\xi_1,\xi_2):=\sup_{y\in Y}\angle_y(\xi_1,\xi_2)$. 
It is useful to know that one can compute $\tangle(\xi_1,\xi_2)$ by only
looking at the angles along a ray $\rho$ asymptotic to one of the
ideal points $\xi_i$;
namely $\angle_{\rho(t)}(\xi_1,\xi_2)$ is monotonically increasing and
converges to $\tangle(\xi_1,\xi_2)$ as $t\to+\infty$.
Another way to represent the Tits metric is as follows:
If $\rho_i$ are rays asymptotic to $\xi_i$, then 
$2\sin\frac{\tangle(\xi_1,\xi_2)}{2}=
\lim_{t\to\infty}\frac{d(\rho_1(t),\rho_2(t))}{t}$. 

The metric space $\tits Y=(\geo Y,\tangle)$ is called 
the {\em Tits boundary}.
It turns out that $\tits Y$ is always a complete metric length space 
with curvature $\leq1$. 
The Tits distance is lower semicontinuous with respect to the cone
topology. 
It induces a topology on $\geo Y$
which is in general strictly finer than the cone topology 
and, generically, $\tits Y$ is not compact
even when $Y$ is locally compact. 
More details can be found in \cite[section 2.3.2]{KleinerLeeb}.

\medskip
Two lines in $Y$ are called {\em parallel} 
if they have finite Hausdorff distance.
Due to a basic rigidity result,
the Flat Strip Lemma, 
any two parallel lines bound an embedded flat strip, that is,
a convex subset isometric to the product of the real line 
with a compact interval. 
The {\em parallel set} $P(l)$ of $l$ is defined as the union of all lines
parallel to $l$.
There is a canonical isometric splitting 
$P(l)\cong l\times CS(l)$
and the {\em cross section} $CS(l)$ is again a Hadamard space.

\medskip
The convexity of the distance $d(\cdot,\cdot)$
provides natural convex functions. 
For instance, 
the distance $d(y,\cdot)$ from a point $y$ is a convex function on $Y$
and, more generally,
the distance $d(C,\cdot)$ from a convex subset $C$.

Related to distance functions are Busemann functions. 
They measure the relative distance from points at infinity.
Their construction goes as follows.
For an ideal point $\xi\in\geo Y$
and a ray $\rho:[0,\8)\to Y$ asymptotic to it
we define the {\em Busemann function} $b_{\xi}$
as the pointwise monotone limit
\[ b_{\xi}(y):=\lim_{t\to\8} (d(y,\rho(t))-t) \]
of normalized distance functions. 
It is a basic but remarkable fact that, up to additive constants, 
$b_{\xi}$ is independent of the ray $\rho$ representing $\xi$. 
As a limit of distance functions, 
$b_{\xi}$ is Lipschitz continuous with Lipschitz constant $1$.
Note that along a ray $\rho$ asymptotic to $\xi$ 
the Busemann function $b_{\xi}$ is affine linear with slope one,
$b_{\xi}(\rho(t)) = -t + const$.

The level and sublevel sets of $b_{\xi}$ are called 
{\em horospheres} and {\em horoballs} centered at $\xi$. 
We denote the horosphere passing through $y$ by $Hs(\xi,y)$ 
and the horoball which it bounds by $Hb(\xi,y)$. 
The horoballs are convex subsets and their ideal boundaries
are convex subsets of $\tits Y$,
namely balls of radius $\pi/2$ around the centers of the horoballs:
$\geo Hb(\xi,y) = 
\{\tangle(\xi,\cdot)\leq\frac{\pi}{2} \}$.

Convex functions have directional derivatives. 
For Busemann functions they are given by the formula
\begin{equation}
\label{busederiv}
\frac{d}{dt^+}(b_{\xi}\circ\si)(t)=-\cos\angle_{\si(t)}(\si'(t),\xi)
\end{equation}
where $\si:I\ra Y$ is a unit speed geodesic segment 
and the angle on the right-hand side is taken 
between the positive direction $\si'(t)\in\Si_{\si(t)}Y$ 
of the segment $\si$ at $\si(t)$ 
and the ray emanating from $\si(t)$ asymptotic to $\xi$.

\subsection{Spherical buildings}
\label{sec:sphbuil}

A {\em spherical Coxeter complex} $(S,W_{sph})$
consists of a unit sphere $S$ and a finite subgroup $W_{sph}\subset Isom(S)$
generated by reflections.
By a reflection, we mean a reflection at a great sphere of codimension one.
$W_{sph}$ is called the {\em Weyl group}
and the fixed point sets of the reflections in $W_{sph}$
are called {\em walls}.
The pattern of walls gives $S$
a natural structure of a cellular (polysimplicial) complex.
The top-dimensional cells, the {\em chambers},
are fundamental domains for the action $W_{sph}\acts S$.
They are spherical simplices if $W_{sph}$ acts without fixed point.
If convenient,
we identify the {\em spherical model Weyl chamber}
$\De_{sph}=S/W_{sph}$
with one of the chambers in $S$.

\medskip
A {\em spherical building} 
modelled on a spherical Coxeter complex $(S,W_{sph})$
is a metric space $B$ with curvature $\leq1$ 
together with a maximal atlas of charts,
i.e.\ isometric embeddings $S\embed B$.
The image of a chart is an {\em apartment} in $B$.
We require that any two points are contained in a common apartment
and that the coordinate changes between charts
are induced by isometries in $W_{sph}$.

We will usually denote the metric on a spherical building by $\tangle$
because in this paper spherical buildings
arise as Tits boundaries of symmetric spaces. 

Two points $\xi,\eta\in B$
are called {\em antipodal}
if $\tangle(\xi,\eta)=\pi$. 

The cell structure and the notions of wall, chamber etc.\
carry over from the Coxeter complex to the building.
The building $B$ is called {\em thick}
if every codimenion-one face is adjacent to at least three chambers.
A non-thick building can always be equipped
with a natural structure of a thick building by reducing the Weyl group.
If $W_{sph}$ acts without fixed points
the chambers are spherical simplices
and the building carries a natural structure
as a piecewise spherical simplicial complex.
We will then refer to the cells as simplices.

There is a canonical 1-Lipschitz continuous
{\em accordion} map $acc:B\ra\De_{sph}$
folding the building onto the model Weyl chamber
so that every chamber projects isometrically.
$acc(\xi)$ is called the {\em type} of the point $\xi\in B$,
and a point in $B$ is called {\em regular}
if its type is an interior point of $\De_{sph}$.

\subsection{Symmetric spaces of noncompact type}
\label{sec:symmsp}

A complete simply connected 
Riemannian manifold $X$ is called a {\em symmetric space}
if in every point $x\in X$ there is a reflection,
that is,
an isometry $\si_x$ fixing $x$ with $d\si_x=-id_x$.
We will always assume that 
$X$ has {\em noncompact type}, 
i.e.\ that it has nonpositive sectional
curvature and no Euclidean factor.
The identity component $G$ of its isometry group 
is then a noncompact semisimple Lie group with trivial center,
and the point stabilizers $K_x$ in $G$ are its maximal compact subgroups.

Given a line $l$ in $X$,
the products of even numbers of reflections at points on $l$
are called {\em translations} or {\em transvections} along $l$.
They form a one-parameter subgroup of isometries in $G$.

With respect to the cone topology
$\ol X$ is a closed standard ball, $X$ its interior
and $\geo X$ the boundary sphere.
The Tits boundary $\tits X=(\geo X,\tangle)$ carries a natural structure 
as a thick {\em spherical building} of dimension $rank(X)-1$.
The faces (simplices) of $\tits X$ correspond to parabolic subgroups of $G$
stabilizing them 
and, as a simplicial complex, 
$\tits X$ is canonically isomorphic 
to the spherical Tits building associated to $G$. 
The building geometry is interesting if $rank(X)\geq2$;
for $rank(X)=1$ the Tits metric 
is discrete with values $0$ and $\pi$.

The top-dimensional simplices of $\tits X$ 
are called {\em (spherical Weyl) chambers}.
They can be simultaneously and compatibly
identified with a {\em spherical model Weyl chamber} $\De_{sph}$. 
In fact,
each orbit for the natural isometric action 
$G\acts\tits X$
meets each chamber in precisely one point
and there is a natural projection  
$acc:\tits X\ra\De_{sph}$
given by dividing out the $G$-action. 
Its restriction to any chamber is an isometry. 
We call $acc$ the {\em accordion map}
because of the way it folds the spherical building 
onto the model chamber.
We refer to the $acc$-image of an ideal point 
as its ($\De_{sph}$-){\em type},
cf.\ section{sec:sphbuil}.

The fixed point set in $\tits X$ of a parabolic subgroup $P\subset G$ 
is a closed simplex $\si_P$.
The stabilizer of each interior point $\xi\in\si_P$
equals $P$
and the map $gP\mapsto g\xi$
defines an embedding 
of the generalized flag manifold $G/P$
into the ideal boundary.
The orbit $G\xi$ is a submanifold of $\geo X$ 
with respect to the cone topology.
If $P$ is a maximal parabolic subgroup
then the fixed point set of $P$ 
is a vertex ($0$-dimensional simplex) of $\tits X$
and we have a unique $G$-equivariant embedding 
$G/P\embed \geo X$.

A subset $Z\subseteq X$ is called a {\em totally-geodesic subspace}
if, with any two distinct points,
it contains the unique line passing through them.
Totally-geodesic subspaces are embedded submanifolds 
and symmetric spaces themselves. 

By a {\em flat} in $X$ we mean a flat totally-geodesic subspace,
i.e.\ a closed convex subset isometric to a Euclidean space.
A $d$-flat is a $d$-dimensional flat.
The Tits metric on $\geo X$ reflects the pattern of flats in $X$:
The Tits boundary 
of a flat $f\subset X$
is a {\em sphere} in $\tits X$, 
by which we mean a closed convex subset 
isometric to a unit sphere in a Euclidean space.
Since $X$ is a symmetric space,
vice versa,
every sphere $s\subset\tits X$ 
arises as the ideal boundary of a flat $f\subset X$,
$\tits f=s$,
actually of several flats if $s$ is not top-dimensional.
The maximal flats in $X$ correspond one-to-one 
to the top-dimensional spheres, the {\em apartments}, in $\tits X$.
The ideal boundary of a {\em singular} flat 
is a subcomplex of $\tits X$. 

The natural action of $G$ 
on the set of all maximal flats is transitive 
and their dimension is called the {\em rank} of the symmetric space.
We have $rank(X)=dim(\tits X)+1$.

A non-maximal flat 
is called {\em singular} if it arises as the intersection of maximal
flats.
Each maximal flat $F$ contains finitely many families of parallel
codimension-one singular flats.
We will also call them {\em singular hyperplanes} in $F$.
Each singular flat $f\subset F$ 
can be obtained as the intersection of singular hyperplanes in $F$.
The ideal boundaries of singular flats in $F$
are subcomplexes of the apartment $\tits F$ in $\tits X$.

For a point $x\in F$,
there are finitely many singular hyperplanes 
$f\subset F$ passing through $x$.
They divide $F$ into cones
whose closures are called {\em Euclidean Weyl chambers}
with {\em tip} $x$.
The reflections at the hyperplanes $f$ generate a finite group, 
the {\em Weyl group} $W_{F,x}$.
It acts on $\tits F$ by isometries 
and $(\tits F,W_{F,x})$ is the {\em spherical Coxeter complex}
attached to $X$ respectively $\tits X$.
It is well-defined up to automorphisms. 

The Euclidean Weyl chambers in $X$ 
can be canonically identified with each other by isometries in $G$, 
and hence they can be simultaneously identified with a 
{\em Euclidean model Weyl chamber} $\De_{euc}$. 
The ideal boundaries of Euclidean Weyl chambers 
are spherical Weyl chambers 
and there is a natural identification 
$\De_{sph}\cong\tits\De_{euc}$.

Let $f$ be a flat, not necessarily singular. 
We call a line $l$ in $f$ {\em maximally regular} 
(with respect to $f$) 
if it is generic in the sense that 
its two ideal endpoints are interior points 
of simplices of $\tits X$
with maximal possible dimension (depending on $f$). 
Every maximal flat $F$ containing $l$ must also contain $f$ 
because the apartment $\tits F$, as a convex subcomplex of $\tits X$, 
must contain the sphere $\tits f$.
Thus the smallest singular flat containing $l$ 
also contains $f$.

\medskip
Any two parallel lines in $X$ bound a flat strip
and, since $X$ is a symmetric space,
are in fact contained in a 2-flat. 
The {\em parallel set} $P(l)$ of a line $l$
is the union of all (maximal) flats containing $l$.
It splits isometrically as 
$P(l)\cong l\times CS(l)$ 
and is a totally-geodesic subspace.
Its cross section $CS(l)$ 
is a symmetric space with $rank(CS(l))=rank(X)-1$. 
The cross section $CS(l)$ has no Euclidean de Rham factor 
if and only if the line is a singular $1$-flat,
equivalently, 
iff its ideal endpoints are vertices of $\tits X$. 

More generally,
we need to consider parallel sets of flats.
Given a flat $f\subset X$,
its parallel set $P(f)$ is defined as the union of all flats $f'$
which are parallel to $f$ in the sense 
that $f$ and $f'$ have finite Hausdorff distance,
or equivalently,
that $\geo f=\geo f'$.
$P(f)$ is the union of all (maximal) flats containing $f$
and it splits isometrically as 
$P(f)\cong f\times CS(f)$.
The cross section $CS(f)$ is a nonpositively curved symmetric space
with $rank(CS(f))=rank(X)-dim(f)$.
It has no Euclidean factor 
if and only if the flat $f$ is singular.
Note that $P(f)$ depends only on the sphere $\tits f$.
If $l$ is a maximally regular line in $f$
then every maximal flat which contains $l$ must also contain $f$
and it follows that $P(f)=P(l)$.

\medskip
The Busemann function $b_{\xi}$ associated to the ideal point $\xi\in\geo X$ 
is smooth.
Its gradient is the unit vector field pointing away from $\xi$,
and the differential is given by 
\begin{equation}
\label{diffbuse}
(db_{\xi})_x(v) =-\cos\angle_x(v,\xi) ,
\end{equation}
compare formula (\ref{busederiv})
in the general case of Hadamard spaces.
The horospheres $Hs(\xi,\cdot)$ centered at $\xi$ are the level sets of
$b_{\xi}$ and thus orthogonal to the geodesics asymptotic to $\xi$.

The Busemann functions are convex, but not strictly convex.
The Hessian $D^2b_{\xi}(x)$ in a point $x$
can be interpreted geometrically as the second fundamental form of the
horosphere $Hs(\xi,x)$. 
The degeneracy of the Hessian is described as follows: 

\begin{lem}
\label{atomicbusemannflinearalonggeo}
Let $u,v\in T_xX$ be non-zero tangent vectors,
let $l_v$ be the geodesic with initial condition $v$,
and suppose that $u$ points towards the ideal point $\xi_u$.
Then the following are equivalent:

(i) 
$D^2_{v,v} b_{\xi_u}=0$.

(ii)
$b_{\xi_u}$ is affine linear on $l_v$.

(iii)
$u$ and $v$ span a $2$-plane in $T_xX$ with sectional curvature zero, 
or they are linearly dependent, 

(iv)
$u$ and $v$ are tangent to a 2-flat or linearly dependent. 
\end{lem}

\begin{ex}
[Busemann functions for the symmetric space associated to $SL(m,\C)$]
\label{ex:busesl}
Let $V$ be a finite-dimensional complex vector space 
equipped with a Hermitian scalar product $((\cdot,\cdot))$. 
Furthermore, 
let $G=SL(V)$, $K$ the maximal compact subgroup preserving 
$((\cdot,\cdot))$ 
and let $X=G/K$ be the associated symmetric space of noncompact type. 

The stabilizer $G_{[v]}\subset G$ of a point $[v]$ in projective space $\P V$
is a maximal parabolic subgroup 
and there is a unique $G$-equivariant embedding 
$\P V\embed\geo X$.
We may hence regard $\P V$ as a $G$-orbit 
(of vertices) in the spherical building $\tits X$.
With respect to the cone topology on $\geo X$,
$\P V$ is an embedded submanifold. 

After suitable normalization (rescaling) of the Riemannian metric on $X$, 
one can express the Busemann function at the ideal boundary point $[v]$ as 
\begin{equation}
\label{busev}
b_{[v]}(gK) :=\log\|g^{-1}v\|
\end{equation}
Note that multiplication of $v$ by a scalar
changes $b_{[v]}$ by an additive constant. 
To justify (\ref{busev}) 
we observe that the right-hand side 
is invariant under the maximal parabolic subgroup $G_v$ 
fixing $[v]$.
Its orbits are the horospheres centered at $[v]$. 
It remains to verify that the right-hand side 
is linear with negative slope 
on some (and hence every) oriented geodesic 
$c:\R\to X$ asymptotic to $[v]$. 
The one-parameter group $(T_t)$ of transvections along $c$
has the following form:
There is a direct sum decomposition 
$V=\langle v\rangle\oplus U$
into common eigenspaces for the $T_t$
such that $T_tv=e^{\la t}$ 
and $T_t|_U=e^{-\la' t}\cdot id_U$
with $\la,\la'>0$ and $\la=dim(U)\cdot\la'$.
Hence 
$\log\|(T_t)^{-1}v\|
=-\la t+const$.
\end{ex}

\subsection{Infinitesimal symmetric spaces}
\label{sec:infsymm}

We keep the notation from the previous section. 
Let $o\in X$ be a base point and $K=Stab_G(o)$ the maximal subgroup fixing it.
The tangent space $T_oX$ 
is canonically identified with the orthogonal complement ${\goth p}$ of ${\goth k}$ in ${\goth g}$ 
with respect to the Killing form:
\[ T_oX\cong{\goth p} \]
$K$ acts on ${\goth p}$ by the restriction of the adjoint
representation. 
It is an orthogonal action with respect to the Riemannian metric 
(and the Killing form).
We denote by $Aff({\goth p})$ the group of transformations on ${\goth p}$ 
generated by $K$ and all translations on ${\goth p}$.
We call the geometry, in the sense of Felix Klein,
consisting of the space ${\goth p}$ and the group $Aff({\goth p})$ 
an {\em infinitesimal symmetric space}.

A {\em flat} in ${\goth p}$
is by definition an affine subspace of the form $z+T_of$
where $f$ is a flat in $X$ passing through $o$
and $z$ is an arbitrary vector in ${\goth p}$.
The flat is called {\em singular} if $f$ is singular.
Singular flats are intersections of maximal flats. 
The maximal flats are of the form $z+{\goth a}$
with ${\goth a}$ a maximal abelian subalgebra contained in ${\goth p}$.

We define the {\em parallel set} of a line $l$ in ${\goth p}$ 
as the union of all (maximal) flats containing $l$.
A parallel set is an affine subspace of the form $z+T_oP(c)$ 
where $z$ is an arbitrary vector
and $P(c)$ the parallel set of a geodesic $c\subset X$ through $o$.

The $K$-orbits in ${\goth p}$ are parametrized by the Euclidean model
Weyl chamber $\De_{euc}$. 
In fact, we can think of $\De_{euc}$ as sitting in an abelian
subalgebra ${\goth a}$ as above,
$\De_{euc}\subset{\goth a}$,  
by identifying it with a Weyl chamber. 
Then each $K$-orbit ${\mathcal O} \subset {\goth p}$ 
meets $\De_{euc}$
in a unique point $h$. 
Due to the natural identifications
$Aff({\goth p})\backslash{\goth p}\times{\goth p}
\cong 
K\backslash{\goth p}
\cong\De_{euc}$ 
we can assign to an oriented geodesic segment $\ol{z_1z_2}$ in ${\goth p}$ 
a vector $\si(z_1,z_2)\in\De_{euc}$
which we call its {\em $\De$-length}. 

The exponential map $exp_o:T_oX\ra X$ yields 
a radial projection
\begin{equation*}
{\goth p}-\{0\}\ra\geo X
\end{equation*}
assigning to a tangent vector $v$ the ideal point 
represented by the geodesic ray with initial condition $v$.
This radial projection restricts 
on the unit sphere $S({\goth p})$ of ${\goth p}$
to a homeomorphism $S({\goth p})\ra\geo X$. 

Note that the infinitesimal symmetric space associated to $X$ 
and all structures which we just defined are,
up to canonical isomorphism,
independent of the base point $o$ and the corresponding splitting 
${\goth g}={\goth k}\oplus{\goth p}$.

\subsection{A transversality result for homogeneous spaces}

In this section,
we provide an auxiliary result of differential-topological nature.

\begin{prop}
\label{transversalityinhomog}
Let $Y$ be a homogeneous space for the Lie group $G$,
and let $Z_1,\dots,Z_n$ be embedded submanifolds.
Then, for almost all $(g_1,\dots,g_n)\in G^n$,
the submanifolds $g_1Z_1,\dots,g_nZ_n$ intersect transversally. 
\end{prop}
\proof
The maps 
$G\times Z_i\lra Y$
are submersions, 
and hence 
the inverse image 
\[ N:=\{(g_1,z_1,\dots,g_n,z_n):g_1z_1=\dots=g_nz_n\} \]
of the small diagonal in $Y^n$ under the canonical map 
$G\times Z_1\times\dots\times G\times Z_n\lra Y^n$
is a submanifold. 
We consider the natural projection 
$N\lra G^n$. 
According to Sard's theorem,
the regular values form a subset of full measure in $G^n$ 
(i.e.\ the set of singular values has zero measure).
Let $(g_1^0,\dots,g_n^0)$ be a regular value. 
Then the intersection $\bigcap g_i^0Z_i$ is a submanifold. 
It remains to verify that the $g_i^0Z_i$ intersect transversally. 

Assume that this is not the case.
It means that there exist $z_i^0\in Z_i$
so that $g_1^0z_1^0=\dots=g_n^0z_n^0=:y$,
and that all $g_i^0Z_i$ are tangent in $y$ to some hypersurface 
$S=\{f=0\}$ with a smooth function $f:Y\ra\R$. 
Define $\psi_i:G\ra\R$ by $\psi_i(g):=f(gz_i^0)$.
Then the composed maps
\[ N\lra Y^n\buildrel{(f,\dots,f)}\over\lra \R^n \]
and 
\[ N\lra G^n\buildrel{(\psi_1,\dots,\psi_n)}\over\lra\R^n \]
have the same differential in the point $p=(g_1^0,z_1^0,\dots,g_n^0,z_n^0)$. 
However,
the image of the first map is at $p$ tangent to the diagonal of $\R^n$, 
while the second map has maximal rank in $p$.
This is a contradiction,
and it follows that the intersection is transversal.
\qed

\begin{rem}
In the algebraic category one can prove a more precise result,
namely that the intersection is transversal 
for a Zariski open subset of tuples $(g_1,\dots,g_n)$,
compare Kleiman's transversality theorem \cite{Kleiman}.
\end{rem}

\section{Stable weighted configurations at infinity}
\label{sec:conf}

We define 
in sections \ref{sec:defstab} and \ref{sec:semi}
a notion of stability for measures and weighted configurations 
on the ideal boundary of a symmetric space $X$ of noncompact type.
This is done 
as in \cite{DouadyEarle}
by associating to 
a measure on $\geo X$ a natural convex function,
a weighted Busemann function.
Stability is then defined in terms of its asymptotic properties. 
As a preparation,
we study in section \ref{sec:asysl} 
properties of convex Lipschitz functions 
on nonpositively curved spaces 
and specialize in section \ref{sec:wbuse}
to weighted Busemann functions on a symmetric space. 
In sections \ref{sec:semistable} and \ref{sec:unstable} 
we investigate properties of measures 
under various stability assumptions.
For instance,
we show the existence of directions of steepest asymptotic descent for
Busemann functions of unstable measures 
and deduce an analogue of the Harder-Narasimhan Lemma
(Theorem \ref{hardernarasimhan}). 
In section \ref{sec:sineq}
we prove Theorem \ref{generalschubertinequalitiesintro}, 
the main result of this paper. 
It provides a finite system of homogeneous linear inequalities 
describing the possible $\De$-weights 
for semistable configurations.

\subsection{Asymptotic slopes of convex functions 
on nonpositively curved spaces} 
\label{sec:asysl}

Let $Y$ be a Hadamard space, 
i.e.\ a space of nonpositive curvature. 
We will now discuss asymptotic properties of Lipschitz continuous 
convex functions $f: Y\ra \R$.
Later on they will be applied to convex combinations of Busemann functions
on symmetric spaces. 

Such a function $f$ is asymptotically linear along any ray
$\rho:[0,\infty)\to Y$.
We define the {\em asymptotic slope} of $f$ at the ideal point
$\eta\in\geo Y$ 
represented by $\rho$ as 
\begin{equation}
\label{slopeconv}
slope_f(\eta):= \lim_{t\ra+\infty} \frac{f(\rho(t))}{t}  .
\end{equation}
That the limit does not depend on the choice of $\rho$
follows, for instance, 
from the Lipschitz assumption. 
Since convex functions of one variable have one-sided derivatives,
we can rewrite (\ref{slopeconv}) as 
\begin{equation*}
slope_f(\eta)= \lim_{t\ra+\infty} \frac{d}{dt^+}(f\circ\rho)(t)
\end{equation*}

\begin{lem}
\label{bdsublev}
For any value $a$ of $f$ holds 
\begin{equation*}
\geo\{f\leq a\} =\{ slope_f\leq 0 \}
\end{equation*}
\end{lem}
\proof
The sublevel set $\{f\leq a\}\subset Y$ is non-empty and convex.
Let $p$ be a point in it. 
For any ideal point $\xi\in\geo\{f\leq a\}$
the ray $\ol{p\xi}$ is contained in $\{f\leq a\}$.
Thus $f$ non-increases along it and $slope_f(\xi)\leq0$.
Vice versa, if $\xi\in\geo Y$ is an ideal point with $slope_f(\xi)\leq0$
then $f$ is non-increasing along $\ol{p\xi}$. 
Hence 
$\ol{p\xi}\subset\{f\leq a\}$ and 
$\xi\in\geo\{f\leq a\}$.
\qed

\medskip
We call 
$\{ slope_f\leq 0 \}$
the set of {\em asymptotic decrease}. 

A subset $C$ of a space with curvature $\leq1$ is called 
{\em convex} if for any two points $p,q\in C$ with $d(p,q)<\pi$ 
the unique shortest segment $\ol{pq}$ is contained in $C$.

\begin{lem}
\label{properconvexfunction}
(i)
The asymptotic slope function 
$slope_f:\tits Y\ra\R$ is Lipschitz continuous 
with the same Lipschitz constant as $f$.

(ii)
The set $\{ slope_f\leq 0 \}\subset\geo Y$ 
is convex with respect to the Tits metric. 
The function $slope_f$ is convex on $\{ slope_f\leq 0 \}$
and strictly convex on $\{ slope_f < 0 \}$.

(iii)
The set $\{ slope_f < 0 \}$ contains no pair of points with distance $\pi$.
If 
it is non-empty
then $slope_f$ has a unique minimum.

(iv)
If $Y$ is locally compact,
then 
$f$ is proper and bounded below if and only if $slope_f>0$ everywhere 
on $\geo Y$. 
\end{lem}
\proof 
(i)
Let $\xi_1,\xi_2\in\tits Y$
and let $\rho_i:[0,+\infty)\ra Y$ be rays asymptotic to $\xi_i$ with the same
initial point $y$.
Then 
$d(\rho_1(t),\rho_2(t))\leq t\cdot 2\sin\frac{\tangle(\xi_1,\xi_2)}{2}
\leq t\cdot \tangle(\xi_1,\xi_2)$.
If $f$ is $L$-Lipschitz, we estimate:
$f(\rho_2(t))\leq f(\rho_1(t))+Lt\cdot\tangle(\xi_1,\xi_2)$,
so 
$\frac{f(\rho_2(t))}{t}\leq
\frac{f(\rho_1(t))}{t}+L\cdot\tangle(\xi_1,\xi_2)$.
Passing to the limit as $t\ra+\infty$ yields the assertion. 

(ii)
Suppose now that $\xi_1,\xi_2\in\geo \{slope_f\leq 0\}$ with 
$\tangle(\xi_1,\xi_2)<\pi$. 
Then the midpoints $m(t)$ of the segments $\ol{\rho_1(t)\rho_2(t)}$ 
converge to the midpoint $\mu$ of $\ol{\xi_1\xi_2}$ in $\tits Y$.
Since $f(m(t))\leq f(y)$,
we have $slope_f(\mu)\leq0$. 
Thus $\{ slope_f\leq 0 \}\subset\geo Y$ is convex. 

In order to estimate the asymptotic slope at $\mu$,
we observe that 
$f(\rho_i(t))\leq slope_f(\xi_i)\cdot t+f(y)$
and thus
\begin{equation}
\label{estimatef}
f(m(t))\leq\frac{slope_f(\xi_1)+slope_f(\xi_2)}{2}\cdot t+f(y) .
\end{equation}
Furthermore
$\lim_{t\to+\infty}\frac{d(\rho_1(t),\rho_2(t))}{t}
=2\sin\frac{\tangle(\xi_1,\xi_2)}{2}$.
The latter fact implies via triangle comparison that 
$\limsup_{t\to+\infty}\frac{d(y,m(t))}{t}
\leq \cos\frac{\tangle(\xi_1,\xi_2)}{2}$.
Using $slope_f(\xi_i)\leq0$ 
we deduce
\begin{equation}
\label{slopeesimate}
slope_f(\mu)\leq
\limsup_{t\to+\infty}\frac{f(m(t))}{d(y,m(t))}\leq 
\frac{slope_f(\xi_1)+slope_f(\xi_2)}{2\cos(\tangle(\xi_1,\xi_2)/2)} ,
\end{equation}
and the convexity properties of $slope_f$ follow.

(iii)
Assume that $\tangle(\xi_1,\xi_2)=\pi$ 
and $slope_f(\xi_i)<0$. 
Then $\frac{d(\rho_1(t),\rho_2(t))}{t}\ra2$
and, by triangle comparison,
$\frac{d(y,m(t))}{t}\ra0$ as $t\to+\infty$.
Since $f$ is Lipschitz, this implies
$\frac{f(m(t))}{t}\ra0$. 
We obtain a contradiction with (\ref{estimatef}),
hence the first assertion holds. 

Suppose that $\eta_n$ are ideal points with 
$slope_f(\eta_n)\ra\inf slope_f<0$. 
Then (\ref{slopeesimate}) implies that the sequence $(\eta_n)$ is Cauchy. 
Since $\tits Y$ is complete, it follows that there is one and only one minimum
for $slope_f$. 

(iv)
If $f$ were not proper,
sublevel sets would be noncompact 
and hence, by local compactness, contain rays.
It follows that there exists an ideal point
with asymptotic slope $\leq0$.
Conversely,
if $\xi$ is an ideal point with $slope_f(\xi)\leq0$
then $f$ is non-increasing along rays asymptotic to $\xi$ and hence not
proper. 
\qed

\medskip
The condition $slope_f\geq0$
does not imply a lower bound for $f$. 
But although there are in general no almost minima
there still exist almost critical points. 
We prove a version in the smooth case. 

\begin{dfn}
A differentiable function $\phi:M\ra\R$ 
on a Riemannian manifold 
is said to have {\em almost critical points}
if there exists a sequence $(p_n)$ of points in $M$
such that $\|\nabla \phi(p_n)\|\to0$.
\end{dfn}

\begin{lem}
\label{semistablealmostcrit}
Let $Y$ be a Hadamard manifold 
and let $f:Y\ra\R$ be a smooth convex function. 
If $slope_f\geq0$
then $f$ has almost critical points. 
\end{lem}
\proof
Suppose that $f$ does not have almost critical points.
This means that there is a lower bound 
$\|\nabla f\|\geq\eps>0$ 
for the length of the gradient of $f$.

We consider the normalized negative gradient flow for $f$,
that is, the flow for the vector field 
$V=-\frac{\nabla f}{\|\nabla f\|}$. 
Its trajectories have unit speed and are complete.
For the derivative of $f$ along a trajectory $\ga:\R\to X$ holds 
\[ (f\circ\ga)'=\langle\nabla f,V\rangle\leq-\eps .\]
We let $y_n:=\ga(n)$ for $n\in\N$ and fix a base point $o\in X$. 
Since $f(y_n)\leq f(o)-n\eps\ra-\infty$,
the points $y_n$ diverge to infinity. 
We connect $o$ to the points $y_n$
by unit speed geodesic segments $\ga_n:[0,l_n]\to X$.
Then $l_n=d(o,y_n)\leq n$.  
Since $X$ is locally compact, 
the sequence of segments $\ga_n$ subconverges to a ray $\rho:[0,+\infty)\ra X$.
Using the convexity of $f$
we obtain for $t\geq0$ the estimate 
$f(\ga_n(t))\leq f(o)-\frac{t}{l_n}n\eps\leq f(o)-t\eps$
and, by passing to the limit,
$f(\rho(t))\leq f(o)-t\eps$.
This implies $slope_f(\eta)\leq-\eps<0$
for the ideal point $\eta$ represented by $\rho$.
\qed

\medskip
The next result compares the asymptotic slopes of convex and linear
functions on flat spaces. 

\begin{lem}
\label{slopeestflat}
Let $E$ be a Euclidean space and 
$f:E\ra\R$ a convex Lipschitz function.
Suppose that $slope_f$ assumes negative values
and let $\xi\in\tits E$ be the unique minimum of $slope_f$.
Then on $\tits E$ holds the inequality:
\begin{equation*}
\label{slopeestimateonflat}
slope_f\geq slope_f(\xi)\cdot\cos\tangle(\xi,\cdot)
\end{equation*}
\end{lem}
\proof
We pick a base point $o$ in $E$
and simplify the function $f$ by a rescaling procedure. 
Consider for $a>0$ the functions $f_a(x):=\frac{1}{a}\cdot f(ax)$ 
where $ax$ denotes the image of $x$ under the homothety with scale factor $a$
and center $o$. 
As $a\ra+\infty$, these functions converge uniformly on compacta to a convex
function $f_{\infty}$ with the same Lipschitz constant as $f$.
Moreover,
$f_{\infty}$ is linear along rays initiating in $o$ and has the same
asymptotic slopes as $f$, i.e.\ 
$slope_{f_{\infty}}\equiv slope_f$ on $\tits E$.
We may assume without loss of generality that $f=f_{\infty}$.

Since $\xi$ is the minimum of $slope_f$, 
we have 
\begin{equation*}
f\geq slope_f(\xi)\cdot d(o,\cdot)
\end{equation*}
with equality along the ray $\rho_{\xi}$ with direction $\xi$ starting in $o$.
Let $\eta\in\tits E$ be another ideal point.
We consider the ray $\rho:[0,+\infty)\ra E$ 
towards $\eta$ initiating in $\rho_{\xi}(t_0)$ for
some $t_0>0$. 
Then 
$f(\rho(t))\geq slope_f(\xi)\cdot d(o,\rho(t))$
for $t\geq0$ with equality in $0$.
Hence we obtain the estimate 
\[
\D_{\dot\rho(0)}f 
\geq 
slope_f(\xi)\cdot \D_{\dot\rho(0)}(d(o,\cdot))
=
slope_f(\xi)\cdot\cos\angle_{\rho(0)}(\xi,\eta)
\]
for the partial derivative of $f$ 
in direction of the unit vector $\dot\rho(0)$. 
Of course, $\angle_{\rho(0)}(\xi,\eta)=\tangle(\xi,\eta)$ because $E$ is
flat. 
The convexity of $f$ implies that 
\[
slope_f(\eta)\geq 
\D_{\dot\rho(0)}f 
\geq
slope_f(\xi)\cdot\cos\tangle(\xi,\eta).
\]
\qed

\subsection{Weighted Busemann functions on symmetric spaces}
\label{sec:wbuse}

Let from now let $X$ denote a symmetric space of noncompact type. 
The class of convex functions on $X$ relevant for this paper 
are finite convex combinations of Busemann functions.
(See section \ref{sec:hadamard}
for the definition of Busemann functions
which we will henceforth also refer to as {\em atomic} Busemann functions.) 
Since it does not complicate the discussion
of their basic properties,
we will consider more generally
continuous convex combinations.

Let ${\cal M}(\geo X)$ 
be the space of Borel measures on $\geo X$ 
with finite total mass 
equipped with the weak $\ast$ topology.
We recall that $\geo X$ carries the cone topology and is homeomorphic to a
sphere. 
The natural $G$-action on ${\cal M}(\geo X)$ is continuous. 
To a measure $\mu\in{\cal M}(\geo X)$ 
we assign the {\em weighted Busemann function} 
\begin{equation*}
b_{\mu} := \int_{\geo X} b_{\xi}\;d\mu(\xi) .
\end{equation*}
It is well-defined up to an additive constant, convex 
and Lipschitz continuous with Lipschitz constant $\|\mu\|$. 

Let $o\in X$ be a base point and 
normalize the Busemann functions by $b_{\xi}(o)=0$. 
Then the map 
$\geo X\times X\ra\R,(\xi,x)\mapsto b_{\xi}(x)$ 
is continuous 
and consequently also the map 
\begin{equation*}
{\cal M}(\geo X)\times X\lra\R,\quad (\mu,x)\mapsto b_{\mu}(x) .
\end{equation*}

\begin{lem}
\label{weightedbusemannflinearalonggeo}
Let $v$ be a non-zero tangent vector 
and $l$ the geodesic with initial condition $v$.
Then the following are equivalent:

(i)
$D^2_{v,v} b_{\mu} =0$.

(ii)
$b_{\mu}$ is affine linear on $l$.

(iii)
$\mu$ is supported on $\geo P_l$. 
\end{lem}
\proof
This follows readily from 
Lemma \ref{atomicbusemannflinearalonggeo},
the corresponding result for atomic Busemann functions, 
because integration yields
$D^2_{v,v} b_{\mu} =\int_{\geo X}D^2_{v,v} b_{\xi}\;d\mu(\xi)$.
Since $D^2_{v,v} b_{\xi}\geq0$,
we have $D^2_{v,v} b_{\mu} =0$ if and only if 
$D^2_{v,v} b_{\xi}=0$ for $\mu$-almost all $\xi$.
Hence (i)$\Ra$(iii)
by Lemma \ref{atomicbusemannflinearalonggeo}.
Clearly (iii)$\Ra$(ii)$\Ra$(i).
\qed

\medskip
We denote by $MIN(\mu)\subset X$ the minimum set of $b_{\mu}$. 
It is convex
but possibly empty.
If $b_{\mu}$ attains a minimum
then $slope_{\mu}\geq0$ everywhere on $\geo X$.
Moreover, 
Lemma \ref{bdsublev} implies that 
\begin{equation*}
\geo MIN(\mu)=\{slope_{\mu}=0\}.
\end{equation*}
Here and later on we abbreviate 
\begin{equation*}
slope_{\mu}:=slope_{b_{\mu}} .
\end{equation*}
By Lemma \ref{weightedbusemannflinearalonggeo}, 
$MIN(\mu)$ contains with any two distinct points 
also the complete geodesic passing through these points. 
Hence:
\begin{cor}
\label{mintotgeod}
If non-empty, $MIN(\mu)$ is a totally geodesic subspace of $X$.
\end{cor}

We compute now the asymptotic slopes of weighted Busemann functions.
A basic observation is that for an atomic Busemann function 
$b_{\xi}:Y\ra\R$ 
the asymptotic slope function on $\tits X$
can be expressed in terms of the Tits geometry.
Using formula (\ref{busederiv}) in section \ref{sec:hadamard} 
for the derivative of Busemann functions 
one obtains 
\begin{equation*}
slope_{b_{\xi}}(\eta)=
\lim_{t\ra+\infty}\frac{d}{dt^+}(b_{\xi}\circ\rho)(t)
=-\lim_{t\ra+\infty}\cos\angle_{\rho(t)}(\xi,\eta)
\end{equation*}
and, 
since $\angle_{\rho(t)}(\xi,\eta)\nearrow\tangle(\xi,\eta)$ as $t\to+\infty$:
\begin{equation*}
slope_{b_{\xi}}(\eta)=-\cos\tangle(\xi,\eta) .
\end{equation*}
The differential of a weighted Busemann function is obtained 
by integrating (\ref{busederiv}): 
\begin{equation*}
(db_{\mu})_x = -\int_{\geo X} \cos\angle_x(\cdot,\xi) \;d\mu(\xi) 
\end{equation*}
With the monotone convergence theorem for integrals we get
\begin{equation}
\label{slopeform}
slope_{\mu} = -\int_{\geo X} \cos\tangle(\cdot,\xi) \;d\mu(\xi).   
\end{equation}
Notice that the asymptotic slope function is expressed
directly in terms of the Tits geometry on the ideal boundary. 

\medskip
We wish to describe the asymptotics of Busemann functions more precisely.
For $\xi,\eta\in\geo X$
and a ray $\rho:[0,+\infty)\ra X$ asymptotic to $\eta$
we have that the convex non-increasing function 
\begin{equation*}
(b_{\xi}\circ\rho)(t) + \cos \tangle(\xi,\eta) \cdot t 
\end{equation*}
converges to a finite limit as $t\to+\infty$.
To see this, we consider 
a flat $F$ with $\xi,\eta\in\geo F$
and inside it a ray $\rho'$ asymptotic to $\eta$. 
Then $b_{\xi}\circ\rho'$ is linear with slope $-\cos\tangle(\xi,\eta)$
and one can estimate 
$|(b_{\xi}\circ\rho)(t)-(b_{\xi}\circ\rho')(t)|\leq d(\rho(t),\rho'(t))
\leq d(\rho(0),\rho'(0))$
because $b_{\xi}$ is 1-Lipschitz.

This kind of asymptotic behavior holds more generally along Weyl chambers. 
Let $F$ be a maximal flat and $V\subset F$ a Weyl chamber.
There exists a maximal flat $F'$ asymptotic to $\xi$ and $V$, i.e.\ with
$\{\xi\}\cup\geo V\subset\geo F'$. 
The restriction of $b_{\xi}$ to $F'$ is then affine linear. 
With a (purely) parabolic isometry $n$ which fixes $\geo V$ and moves $F$ to
$F'$, we may write  
\begin{equation}
\label{decomposeb}
b_{\xi}(x)=b_{\xi}(nx)+(b_{\xi}(x)-b_{\xi}(nx)) .
\end{equation}
The summand $b_{\xi}(nx)$ is linear on $V$ whereas $b_{\xi}(x)-b_{\xi}(nx)$
is bounded (and convex) on $V$. 
Thus the restriction of a Busemann function $b_{\xi}$ to a Euclidean Weyl
chamber is {\em asymptotically linear} in the sense that it is the sum of a
linear and a bounded function. 

We generalize to weighted Busemann functions:

\begin{lem}[Asymptotic linearity]
\label{aslin}
The Busemann function $b_{\mu}$ is asymptotically linear on each Euclidean Weyl
chamber $V\subset X$ 
in the sense that its restriction to $V$ decomposes as the sum of a linear
function and a sublinear function. 
If the support of $\mu$ is finite, the sublinear part is bounded.
\end{lem}
\proof 
According to (\ref{decomposeb})
we can decompose each atomic Busemann function $b_{\xi}$ on $V$ as the sum 
$b_{\xi}|_V = l_{\xi} + s_{\xi}$
of its linear and bounded part. 
For measures $\mu$ with finite support the claim follows directly. 

For arbitrary measures $\mu\in{\cal M}(\geo X)$ one has to argue a bit more
carefully. 
We normalize the functions $b_{\xi}$, $l_{\xi}$ and $s_{\xi}$ 
to be zero at the tip of $V$. 
The decomposition 
of $b_{\xi}$ depends measurably on $\xi$.
Note that $l_{\xi}$ is 1-Lipschitz and hence $s_{\xi}$ is 2-Lipschitz. 
This allows us to integrate 
and we get 
$b_{\mu}=\int l_{\xi}\;d\mu(\xi)+\int s_{\xi}\;d\mu(\xi)$.
Both summands are Lipschitz and the first one is clearly linear.
To see that the second one is sublinear 
one uses that $s_{\xi}$ is bounded for each $\xi$ 
and $\sup_V s_{\xi}$ depends measurably on $\xi$. 
\qed

\medskip
As a consequence of asymptotic linearity, 
the function $slope_{\mu}$ has the property
that its values on a simplex $\si\subset\tits X$
are determined by the values on the vertices of $\si$.
Namely, if $F$ is a flat in $X$ with $\si\subset\tits F$
and if $l:F\ra\R$ is an affine linear function 
with the same asymptotic slopes at the vertices of $\si$ as $b_{\mu}$
then $slope_{\mu}=slope_l$ on $\si$.
Since $\{slope_l>0\}$ resp.\ $\{slope_l\geq0\}$
is an open resp.\ closed hemisphere in the round sphere $\tits F$,
this implies:

\begin{cor} 
\label{vertexslopes}
(i) 
Suppose that $slope_{\mu} > 0$ 
(resp. $slope_{\mu} \geq 0$, 
$slope_{\mu} \leq 0$ or
$slope_{\mu} <0$) 
on all vertices of a simplex $\si\subset\tits X$. 
Then the same inequality holds on the entire simplex $\si$.

(ii) If $slope_{\mu} \geq 0$
holds on $\si$
then $\{slope_{\mu}=0\}\cap\si$ is a face of $\si$.
\end{cor}

\subsection{Stability for measures on the ideal boundary}
\label{sec:defstab}

We define stability of a measure $\mu\in{\cal M}(\geo X)$ 
in terms of its weighted Busemann function $b_{\mu}$ on $X$ 
and the associated asymptotic slope function on $\tits X$. 

\begin{dfn}[(Semi)Stability of measures on $\tits X$]
\label{def:stabmeasinf}
We call a measure $\mu\in{\cal M}(\geo X)$ 
{\em stable} if $slope_{\mu}>0$, 
{\em semistable} if $slope_{\mu}\geq0$ 
and 
{\em unstable} if it is not semistable.
\end{dfn}

\begin{rem}
In fact,
formula (\ref{slopeform})
expresses $slope_{\mu}$ directly in terms of the 
intrinsic geometry of $\tits X$ 
without referring to $b_{\mu}$.
Our definition of stability hence 
carries over to Borel measures with finite total mass 
on topological spherical buildings 
in the sense of Burns and Spatzier \cite{BurnsSpatzier}. 
It agrees with \cite[Definition 4.1]{polbuil}
given in the special case of measures with finite support. 
In this case the integration in (\ref{slopeform})
becomes finite summation 
and makes sense on any spherical building
(which may be thought of as a topological spherical building
with discrete topology). 
\end{rem}

If the measure $\mu$ is semistable 
then $\{slope_{\mu}=0\}$ 
is a convex subcomplex of $\tits X$ 
by Lemma \ref{properconvexfunction} (ii) 
and Corollary \ref{vertexslopes} (ii).
The following more subtle variation 
of the notion of stability 
will be needed in section \ref{sec:polyxconf}, 
in particular for Proposition \ref{prop:stabletofixed}
and the proof of Theorem \ref{polxconf:samelengths}. 

\begin{dfn}[Nice semistability of measures on $\tits X$]
\label{def:nss}
We call a semistable measure $\mu$ 
{\em nice semistable} if 
$\{slope_{\mu}=0\}$
is either empty or $d$-dimensional and contains a unit $d$-sphere. 
\end{dfn}

\begin{rem}
A $d$-dimensional convex subcomplex of a spherical building 
which contains a unit $d$-sphere
carries itself a natural structure as a spherical building. 
In fact we will show in Lemma \ref{nss} that 
for a nice semistable measure $\mu$ 
the set $\{slope_{\mu}=0\}$
is the ideal boundary of a totally geodesic subspace.
\end{rem}

In view of the slope formula (\ref{slopeform})
semistability of $\mu$ 
is equivalent to the system of inequalities 
\begin{equation}
\label{ssineq}
\int_{\geo X} \cos\tangle(\eta,\cdot) \;d\mu  \leq 0
\qquad\qquad\forall\;\eta\in\geo X 
\end{equation}
and stability to the corresponding system of strict inequalities. 
Note that according to asymptotic linearity 
(Corollary \ref{vertexslopes})
it suffices to check the inequalities on vertices.

\begin{ex}
\label{ex:stabrank1}
Suppose that $X$ has rank one,
equivalently, 
that the spherical building $\tits X$ has dimension $0$.
Then a measure $\mu$ on $\tits X$ is stable
if and only if it has no atom with mass $\geq\half|\mu|$,
semistable 
if and only if it has no atom with mass $>\half|\mu|$,
and nice semistable if and only if 
it is either stable or consists of two atoms with equal mass. 
This follows from the fact that the Tits metric is discrete 
(with distances $0$ or $\pi$) 
and hence
\begin{equation*}
slope_{\mu}(\eta)
=
\int_{\geo X} \cos\tangle(\eta,\cdot) \;d\mu 
=
1\cdot\mu(\eta)+(-1)\cdot(|\mu|-\mu(\eta))
=
2\cdot\mu(\eta)-|\mu|.
\end{equation*}
\end{ex}
In higher rank 
the stability criterion becomes more complicated.
Examples are discussed in section \ref{sec:expl}
where we work out the case of measures on the 
Grassmannians associated to the classical groups.

\medskip
The next result implies 
that semistability persists under totally
geodesic embeddings of symmetric spaces. 
It is useful when relating the stability conditions for different groups.

\begin{lem}
\label{masssuportedonsubspace}
Assume that $C\subset X$ is a closed convex subset 
and that the measure $\mu$ is supported on $\geo C\subset\geo X$.
Then:

(i)
$\inf b_{\mu}\restr_C=\inf b_{\mu}$.
In particular, 
$b_{\mu}$ is bounded below on $C$
if and only if it is bounded below on $X$.

(ii)
If $slope_{\mu}\geq0$ on $\geo C$
then $slope_{\mu}\geq0$ on $\geo X$.
\end{lem}
\proof
For every ideal point $\xi\in\geo C$ holds 
$b_{\xi}\geq b_{\xi}\circ\pi_C$ 
where $\pi_C:X\ra C$ 
denotes the nearest point projection.
Namely, let $x$ be a point in $X$ 
and let $\si$ be the segment connecting $x$ and its projection $\pi_C(x)$.
Then for any point $x'$ on $\si$
we have that 
the ideal triangle with vertices $x',\pi_C(x)$ and
$\xi$ has angle $\frac{\pi}{2}$ at $\pi_C(x)$ and therefore angle
$\leq\frac{\pi}{2}$ at $x'$ 
because the angle sum is $\leq\pi$. 
Hence $b_{\xi}$ decreases along $\si$ and we obtain 
$b_{\xi}(x)\geq b_{\xi}(\pi_C(x))$. 
Integration with respect to $\mu$ yields:
\begin{equation}
\label{buseproj}
b_{\mu}\geq b_{\mu}\circ\pi_C 
\end{equation}
This implies assertion (i).

Regarding part (ii), 
suppose that $slope_{\mu}\geq0$ on $\geo C$ and $slope_{\mu}(\eta)<0$ for
some $\eta\in\geo X$.
Let $\rho:[0,+\infty)\ra X$ be a unit speed ray asymptotic to $\eta$
with $b_{\mu}(\rho(0))\leq0$.
Then $b_{\mu}(\rho(t))\leq -ct$ with $c:=-slope(\eta)>0$.
Let $y_n:=\pi_C(\rho(n))$ for $n\in\N_0$.
In view of (\ref{buseproj})
we have $b_{\mu}(y_n)\leq b_{\mu}(\rho(n))\leq -cn$.

Nearest point projections to closed convex sets are $1$-Lipschitz
and therefore $d(y_0,y_n)\leq n$. 
On the other hand $d(y_0,y_n)\ra\infty$ because 
$b_{\mu}(y_n)\ra-\infty$. 
Thus the sequence of segments $\ol{y_0y_n}$ in $C$ subconverges to a ray
$\bar\rho$ in $C$. 
Using the convexity of Busemann functions, 
it follows that 
$b_{\mu}(\bar\rho(t))\leq b_{\mu}(y_0)-ct$
and hence $slope_{\mu}(\bar\eta)\leq-c<0$ at the ideal endpoint
$\bar\eta$ of $\bar\rho$. 
This is a contradiction because $\bar\eta\in\geo C$.
\qed

\subsection{Properties of stable and semistable measures}
\label{sec:semistable}

We investigate now
how the various degrees of stability of a measure 
are reflected in the behavior 
of the associated weighted Busemann function.

\begin{lem}
\label{lem:stable}
$\mu$ is stable
if and only if $b_{\mu}$ is proper and bounded below.
In this case $b_{\mu}$ has a unique minimum. 
\end{lem}
\proof
Part (iv) of Lemma \ref{properconvexfunction} implies 
that $\mu$ is stable if and only if $b_{\mu}$ is proper and bounded below. 
Since $b_{\mu}$ is convex this is in turn equivalent 
to $MIN(\mu)$ being compact and non-empty,
and by Corollary \ref{mintotgeod} 
to $MIN(\mu)$ being a point. 
\qed

\begin{lem}
\label{nss}
$\mu$ is nice semistable
if and only if $b_{\mu}$ attains a minimum. 
\end{lem}
\proof
``$\Rightarrow$'':
Suppose that $\mu$ is nice semistable.
We are done by Lemma \ref{lem:stable} if $\mu$ is stable.
Therefore we assume also that $\{slope_{\mu}=0\}$ 
is non-empty 
and hence a $d$-dimensional convex subcomplex 
which contains a unit $d$-sphere $s$.
Let $f\subset X$ be a flat with $\geo f=s$.
Furthermore,
let $l$ be a maximally regular geodesic inside $f$.
Then any geodesic parallel to $l$ lies in a flat parallel to $f$
and the parallel sets satisfy
$P(f)=P(l)$.
Lemma \ref{weightedbusemannflinearalonggeo} 
implies that $\mu$ is supported on $\geo P(f)$. 
By Lemma \ref{masssuportedonsubspace} 
it suffices to show that the restriction of 
$b_{\mu}$ to the parallel set $P(f)\cong f\times CS(f)$
attains a minimum. 
Since $b_{\mu}$ is constant on each flat parallel to $f$
this amounts to finding a minimum on a cross section 
$\{pt\}\times CS(f)$. 
We have $slope_{\mu}>0$ on $\geo (\{pt\}\times CS(f))$ 
because otherwise $\{slope_{\mu}=0\}$ 
would contain a $(d+1)$-dimensional hemisphere,
which is absurd. 
Using Lemma \ref{properconvexfunction} (iv)
we conclude that $b_{\mu}$ attains a minimum on $\{pt\}\times CS(f)$. 

``$\Leftarrow$'':
If $b_{\mu}$ attains a minimum 
then $\mu$ is semistable 
and $\{slope_{\mu}=0\}=\geo MIN(\mu)$ 
is empty or the ideal boundary of a totally geodesic subspace,
cf.\ Corollary \ref{mintotgeod}. 
In the latter case
$\{slope_{\mu}=0\}$ carries a natural structure as a spherical building 
and hence contains a top-dimensional unit sphere. 
\qed

\medskip
As a special case of Lemma \ref{semistablealmostcrit} we obtain: 

\begin{lem}
\label{almostcritical}
If $\mu$ is semistable 
then $b_{\mu}$ has almost critical points. 
\end{lem}

\begin{lem}
\label{nssinclosure}
If $\mu$ is semistable
then the closure of its $G$-orbit in ${\cal M}(\geo X)$ 
contains a nice semistable measure. 
\end{lem}
\proof
By Lemma \ref{almostcritical}
the associated weighted Busemann function $b_{\mu}$
has almost critical points,
i.e.\ there exists a sequence $(x_j)$ of points in $X$
with 
$\|\nabla b_{\mu}(x_j)\|\to0$. 
We use the $G$-action on $X$ 
to move the almost critical points into the base point. 
Namely let $g_j\in G$ with $g_jx_j=o$.
Then 
$\|\nabla b_{g_j\mu}(o)\|=\|\nabla b_{\mu}(x_j)\|\to0$. 
Due to the compactness of ${\cal M}(\geo X)$
we may assume after passing to a subsequence 
that $g_j\mu\ra\nu$.
It follows that $\nabla b_{g_j\mu}(o)\ra \nabla b_{\nu}(o)$
and therefore $\nabla b_{\nu}(o)$.
Hence $\nu\in \ol{G\mu}$ 
is a nice semistable measure.
\qed

\begin{rem}
(i)
One can show that the closure $\ol{G\mu}$ of a semistable orbit
contains a {\em unique} nice semistable $G$-orbit $G\nu$.

(ii)
The Busemann function of a semistable measure $\mu$ 
is in general not bounded below.
However, 
for semistable measures $\mu$ with {\em finite support}
one can show that $b_{\mu}$ is bounded below, 
but this fact will not be needed in this paper. 
\end{rem}

\subsection{Unstable measures and directions of steepest descent}
\label{sec:unstable}

Lemma \ref{properconvexfunction} implies 
that for unstable measures $\mu$ 
there is a {\em unique} ideal point $\xi_{min}$ of {\em steepest descent},
i.e.\ where $slope_{\mu}$ attains its minimum. 
We will now look for {\em vertices} of steepest descent among vertices of a
given type.
The following uniqueness result is a version of the Harder-Narasimhan Lemma. 

\begin{thm}
\label{hardernarasimhan}
Let $\mu\in{\cal M}(\geo X)$ be unstable,
and let $\xi_{min}$ be the unique ideal point of steepest descent for
$b_{\mu}$. 
Let $\tau_{min}$ be the simplex in the Tits boundary spanned by
$\xi_{min}$, 
i.e.\ which contains $\xi_{min}$ as an interior point. 

Then for each vertex $\eta$ of $\tau_{min}$ holds:
$\eta$ is the unique minimum of $slope_{\mu}$ restricted to the orbit $G\eta$,
i.e.\ the unique minimum among vertices of the same type. 
\end{thm}
Note that $slope_{\mu}<0$ on all vertices of $\tau_{min}$ 
due to the asymptotic linearity of Busemann functions on Weyl chambers 
(Lemma \ref{aslin})
and the fact that all simplices in the Tits boundary have diameter $\leq\frac{\pi}{2}$. 
\proof
Since $slope_{\mu}(\xi_{min})<0$, 
there is a unique measure $\nu$ supported on the vertices of $\tau_{min}$ such
that $\xi_{min}$ is the ideal point of steepest $\nu$-descent and 
$slope_{\nu}(\xi_{min})=slope_{\mu}(\xi_{min})$.
On each flat $f$ asymptotic to $\tau_{min}$ 
i.e.\ with $\tau_{min}\subset\geo f$,
the Busemann function $b_{\nu}$ restricts to a linear function 
whose negative gradient points towards $\xi_{min}$. 
The asymptotic slopes are given by the formula:
\begin{equation*}
slope_{\nu} = slope_{\mu}(\xi_{min})\cdot\cos\tangle(\xi_{min},\cdot) 
\qquad\hbox{on $\tits X$}
\end{equation*}
Let $f$ be a minimal flat containing $\tau_{min}$ in its ideal boundary,
i.e.\ $dim(f)=dim(\tau_{min})+1$ and $\tits f$ is a subcomplex of $\tits X$ 
with $\tau_{min}$ as a top-dimensional simplex.
From the asymptotic linearity of Busemann functions on Weyl chambers
(Lemma \ref{aslin}) 
follows the existence of a linear function $l$ on $f$ 
with $slope_l=slope_{\mu}$ on $\tau_{min}$. 
Since $\xi_{min}$ is the direction of steepest descent for $l$, 
we have that 
$l=b_{\nu}|_f$
modulo additive constants. 
Thus: 
\[ 
slope_{\mu}=slope_{\nu} 
\qquad\qquad\hbox{on $\tau_{min}$}
\]
Every ideal point lies in an apartment through $\xi_{min}$.
Therefore 
by applying Lemma \ref{slopeestflat} to all flats which contain $\xi_{min}$
in their ideal boundary, 
we obtain the estimate:
\begin{equation*}
\label{munusleverywh}
slope_{\mu}\geq slope_{\nu} 
\qquad\qquad\hbox{on $\tits X$}
\end{equation*}
As a consequence,
it suffices to prove: 
{\em (*) $\eta$ is the unique vertex in the orbit $G\eta$ with minimal Tits
distance from $\xi_{min}$.}

To verify this claim,
consider a vertex $\zeta\in G\eta$ in the same orbit.
There exists an apartment $a$ in $\tits X$ containing $\zeta$ and
$\tau_{min}$. 
Suppose that 
$\zeta$ is separated from $\xi_{min}$ inside $a$ by a wall $s$.
The reflection at $s$ belongs to the Weyl group $W(a)$,
and the mirror image $\zeta'$ of $\zeta$ is a vertex of the same type
which is strictly closer to $\xi_{min}$. 
Observe that the vertices which cannot be separated from $\xi_{min}$ by a wall
are precisely the vertices of Weyl chambers with
$\tau_{min}$ as a face. 
Therefore $\eta$ is the only vertex in $G\eta\cap a$ which cannot be
separated from $\xi_{min}$. 
Hence $\eta$ is closer to $\xi_{min}$ than any other vertex in 
$G\eta\cap a$. 
This shows (*) and finishes the proof of Theorem \ref{hardernarasimhan}.
\qed

\medskip
Depending on the geometry of the spherical Weyl chamber and the type of
$\xi_{min}$ one can say more.
See section \ref{sec:symmsp}
for the definition of 
the accordion map $acc:\tits X\ra\De_{sph}$ 
and the definition of type.

\begin{add}
\label{add:steepestvertex}
Suppose that the vertices of $\De_{sph}$ closest to the type
$acc(\xi_{min})$ of $\xi_{min}$ 
belong to the face $acc(\tau_{min})$ spanned by it. 

Then there exists a vertex with steepest $\mu$-descent among all vertices.
All such vertices are vertices of $\tau_{min}$. 
In particular, there are only finitely many of them and each $G$-orbit
contains at most one.
\end{add}
\proof
We take up our argument for Theorem \ref{hardernarasimhan}. 
As we saw,
the vertices at minimal Tits distance from $\xi_{min}$ belong to a Weyl chamber
containing $\tau_{min}$ as a face. 
By our assumption, they are vertices of $\tau_{min}$. 
It follows that the vertices of $\tau_{min}$ closest to $\xi_{min}$ have
minimal $\mu$-slope among vertices. 
\qed

\medskip
Other than one might first expect,
the assumption of \ref{add:steepestvertex} does not always hold.
This can occur if the Dynkin diagram associated to $G$ branches, 
i.e.\ (in the irreducible case) the associated root system is of type 
$D_n$, $E_6$, $E_7$ or $E_8$. 
We discuss the simplest example.

\begin{ex}
Suppose that the Dynkin diagram associated to $G$ is $D_4$.
Then $X$ has rank 4 and the spherical model Weyl chamber $\De_{sph}$ 
is a three-dimensional spherical tetrahedron $\eta\xi_1\xi_2\xi_3$ 
with the following geometry: 
The face $\xi_1\xi_2\xi_3$, which we regard as the base of the tetrahedron, 
is an equilateral triangle.
The dihedral angles at the edges of the base triangle equal $\frac{\pi}{3}$
whereas the dihedral angles at the other three edges $\ol{\eta\xi_i}$ equal
$\frac{\pi}{2}$. 
From these data one deduces that the edges $\ol{\xi_i\xi_j}$, $i\neq j$, have
length $\frac{\pi}{3}$ and the height, i.e.\ the distance between the vertex
$\eta$ and the center $\zeta$ of the base, equals $\frac{\pi}{6}$. 
This shows that $\zeta$ is strictly closer to the opposite vertex $\eta$ than
to the vertices $\xi_i$ of the base.

Consider an atomic measure $\mu$ 
with mass concentrated in one
ideal point $\hat\zeta\in\geo X$ of type $\zeta$.
There are infinitely many vertices of steepest $\mu$-descent among vertices,
namely all vertices $\hat\eta$ 
which are vertices of a Weyl chamber $\si$ such that $\hat\zeta$ is the center
of the face of $\si$ opposite to $\hat\eta$. 

By slightly modifying the measure $\mu$ we can construct a measure $\mu'$ with
infinite support which has no vertices of steepest descent among vertices at
all.
Namely let $A$ be the set consisting of all points $\xi\in\tits X$ with 
$\tangle(\hat\zeta,\xi)=\frac{\pi}{2}$ and such that 
the segment $\ol{\hat\zeta\xi}$ is perpendicular to the simplex $\tau$ spanned
by $\hat\zeta$. 
Then add to $\mu$ a small amount of mass which is distributed over infinitely
many points of $A$. 
\end{ex}

\subsection{Weighted configurations on $\tits X$ and stability}
\label{sec:semi}

A collection 
of points $\xi_1,\dots,\xi_n\in \tits X$
and of weights $m_1,\dots,m_n\geq0$
determines a {\em weighted configuration}
\begin{equation*}
\psi: (\Z/n\Z, \nu)\to \tits X
\end{equation*}
on $\tits X$.
Here $\nu$ is the measure on $\Z/n\Z$ defined by $\nu(i)=m_i$,
and $\psi(i)=\xi_i$. 
By composing $\psi$ with $acc:\tits X\ra\De_{sph}$ 
one obtains a map $(\Z/n\Z,\nu)\ra\De_{sph}$.
It corresponds to a point 
$h=(h_1,\dots,h_n)$
in $\De_{euc}^n$ 
which we call the $\De$-{\em weights} of the configuration $\psi$,
i.e.\ $h_i=m_i\cdot acc(\xi_i)$. 

The configuration $\psi$ yields, by pushing forward $\nu$, 
the measure $\mu=\sum m_i\de_{\xi_i}$ on $\tits X$. 
Accordingly, Definition \ref{def:stabmeasinf} 
carries over from measures to configurations:

\begin{dfn}
[Stability of weighted configurations on $\tits X$]
\label{def:stabconfbuil}
The weighted configuration $\psi$ is called 
stable, semistable, unstable resp.\ nice semistable
if the associated measure $\mu$ has this property.
\end{dfn}

\begin{rem}
Obviously,
the definition extends to weighted configurations 
on abstract spherical buildings.
One may extend it further 
to weighted configurations with infinite support
$\psi:(\Om,\nu)\to B$ on topological spherical buildings $B$,
for instance, on $\tits X$. 
Here $(\Om,\nu)$ denotes a measure space with finite total mass
and the map $\psi$ is supposed to be measurable
with respect to the Borel $\si$-algebra on $B$.
\end{rem}

This notion of stability is motivated by Mumford stability in 
geometric invariant theory. 
The connection between the two concepts is explained in section \ref{git}.

\subsection{The stability inequalities for $\De$-weights of configurations}
\label{sec:sineq}

We will now address the question 
which $\De$-weights occur for semistable weighted configurations 
on $\tits X$. 
We will need the {\em Schubert calculus}.

Think of the model spherical Weyl chamber $\De_{sph}$ as being embedded 
in the spherical Coxeter complex $(S,W)$. 
For a vertex $\zeta$ of $\De_{sph}$,
we denote by $Grass_{\zeta}$ the corresponding maximally singular $G$-orbit 
in $\geo X$.
We call it a {\em generalized Grassmannian}
because in the case of $SL(n)$ the $Grass_{\zeta}$ are the Grassmann manifolds. 
The action of a Borel subgroup $B\subset G$ stratifies each $Grass_{\zeta}$ 
into {\em Schubert cells},
one cell $C_{\eta_i}$ corresponding to each vertex $\eta_i\in S$ 
in the orbit $W\zeta$ of $\zeta$ under the Weyl group $W$. 
Hence, if we denote $W_{\zeta}:=Stab_W(\zeta)$ 
then the Schubert cells correspond to cosets in $W/W_{\zeta}$. 
The {\em Schubert cycles} are defined as the closures $\ol C_{\eta_i}$ 
of the Schubert cells;
they are unions of Schubert cells.
There is one top-dimensional Schubert cell 
corresponding to the vertex in $S$ belonging to the chamber opposite to $\De_{sph}$. 
Note that as real algebraic varieties, 
the Schubert cycles represent homology classes 
$[\ol C_{\eta_i}]\in H^{\ast}(Grass_{\zeta};\Z/2\Z)$
which we abbreviate to $[C_{\eta_i}]$. 
In the complex case they even represent {\em integral} homology classes. 

It will be useful to have another description of the Schubert cells and
Schubert cycles. We recall the definition of the 
{\em relative position} of
a spherical Weyl chamber $\sigma$ and a vertex $\eta$ of 
$\tits X$.  
There exists an apartment $a$ in $\tits X$ containing $\sigma$ and $\eta$.
Furthermore there exists a unique apartment chart $\phi:a\ra S$ which maps
$\sigma$ to $\De_{sph}$.  
We then define the relative position 
$(\sigma,\eta)$ to be the vertex $\phi(\eta)$ of the model apartment $S$. 
To see that the relative position is well-defined, we choose an interior point
$\xi$ in $\si$
and a minimizing geodesic $\ol{\xi\eta}$.
(It is unique if $\tangle(\xi,\eta)<\pi$.)
We then observe that the $\phi$-image 
of the geodesic $\ol{\xi\eta}$ 
is determined by its length and its initial direction in
$\phi(\xi)$,
because geodesics in the unit sphere $S$ do not branch. 
Thus its endpoint $\phi(\eta)$ is uniquely determined by $\si$ and $\eta$.

Notice that $G$ acts transitively on pairs $(\si,\eta)$ with the same relative
position. 
This follows from the transitivity of the $G$-action on pairs $(\si,a)$ of
chambers and apartments containing them.
This implies that the relative position
determines the Tits distance:
\begin{lem} 
\label{relpostitsdis}
Suppose that $\si_1,\eta_1$ and $\si_2,\eta_2$ are given with 
$(\si_1,\eta_1)=(\si_2,\eta_2)$. Suppose further we are given
$\xi_1\in \si_1$ and $\xi_2\in \si_2$ with
$acc(\xi_1)=acc(\xi_2)$. 
Then
\[\tangle(\xi_1,\eta_1)=\tangle(\xi_2,\eta_2).\]
\end{lem}
We now have another description of the Schubert cells as mentioned in the
introduction. As above we assume we have chosen a spherical Weyl chamber
$\si\subset\tits X$ 
and a vertex $\zeta$ of $\De_{sph}$. 
For $\eta_i\in W\zeta$ we then have:
\begin{lem}
The Schubert cell $C_{\eta_i}$ is given by
\[C_{\eta_i}=\{\eta\in Grass_{\zeta}: (\si,\eta)=\eta_i\}.\]
 \end{lem}
For all vertices $\zeta$ of $\De_{sph}$ 
and all $n$-tuples of vertices $\eta_1,\dots,\eta_n\in W\zeta$ 
we consider the inequality
\begin{equation}
\label{stabineq}
\sum_i m_i\cdot\cos\angle(\tau_i,\eta_i)\leq0 ,
\end{equation}
for $m_i\in\R^+_0$ and $\tau_i\in\De_{sph}$
where $\angle$ measures the spherical distance in $S$.
We may rewrite the inequality as follows 
using standard terminology of Lie theory: 
Let $\la_{\zeta}\in\De_{euc}$ be the fundamental coweight 
contained in the edge with direction $\zeta$,
and let $\la_i:=w_i\la_{\zeta}$ where $[w_i]\in W/W_{\zeta}$ 
such that $w_i\zeta=\eta_i$. 
With the renaming $h_i=m_i\tau_i$ of the variables 
(\ref{stabineq}) becomes the homogeneous linear inequality
\begin{equation}
\label{stabineqlie}
\sum_i \langle h_i,\la_i\rangle \leq0 .
\end{equation}
We will now prove our main result 
Theorem \ref{generalschubertinequalitiesintro}
which describes,
in terms of the Schubert calculus, 
a subset of these inequalities 
which is equivalent to the existence of a semistable weighted configurations 
for the given $\De$-weights. 
Theorem \ref{generalschubertinequalitiesintro}
is the combination of the next two theorems. 

\begin{thm}
[Stability inequalities for noncompact semisimple Lie groups]
\label{generalschubertinequalities}
For $h\in\De_{euc}^n$ 
there exists a semistable weighted configuration with $\De$-weights $h$
if and only if 
(\ref{stabineqlie})
holds whenever the intersection 
of the Schubert classes 
$[C_{\eta_1}],\dots,[C_{\eta_n}]$ 
in $H_{\ast}(Grass_{\zeta};\Z/2\Z)$ 
equals $[pt]$.
\end{thm}
\proof
"$\Leftarrow$":
Assume that all configurations with $\De$-weights $h$ are unstable.
Due to the transversality result \ref{transversalityinhomog},
there exist chambers 
$\si_1,\dots,\si_n\subset\tits X$ so that 
the corresponding $n$ stratifications of the Grassmannians $Grass_{\zeta}$ 
by orbits of the Borel subgroups $B_i=Stab_G(\si_i)$ are transversal. 
(This transversality is actually generic.) 
We choose a configuration with $\De$-weights $h$ so that the atoms $\xi_i$ 
are located on the chambers $\si_i$.

Now we apply the Harder-Narasimhan Lemma type 
Theorem \ref{hardernarasimhan}. 
Since the measure $\mu$ on $\geo X$ associated to the configuration is unstable
there exists a vertex $\zeta$ of $\De_{sph}$ 
such that on the corresponding Grassmannian $Grass_{\zeta}$ there is a 
{\em unique} minimum $\eta_{sing}$ for $slope_{\mu}$. 

Let $C_i:=B_i\cdot\eta_{sing}$ be the Schubert cell 
passing through $\eta_{sing}$ 
for the stratification of $Grass_{\zeta}$ by $B_i$-orbits. 
Note that 
all points in the intersection $C_1\cap\dots\cap C_n$ have the same relative
position with respect to all atoms $\xi_i$ and therefore they have equal
$\mu$-slopes. 
Since $\eta_{sing}$ is the unique minimum of $slope_{\mu}$ on $Grass_{\zeta}$ 
it is hence the unique intersection point of the Schubert cells $C_i$. 

Transversality implies that the corresponding Schubert cycles $\bar C_i$
intersect transversally in the unique point $\eta_{sing}$.
The corresponding inequality 
(\ref{stabineq}) resp.\ 
(\ref{stabineqlie})
in our list is violated
because the left sides equal $-slope_{\mu}(\eta_{sing})>0$.

"$\Ra$":
Conversely,
assume that there exists a semistable configuration $\psi$ on $\tits X$
with $\De$-weights $h$
and masses $m_i=\|h_i\|$ located in the ideal points $\xi_i$
of type $\tau_i=\frac{h_i}{\|h_i\|}=acc(\xi)$.
Assume further that 
we have a homologically non-trivial product of $n$ Schubert classes:
$[C_{\eta_1}]\cdots[C_{\eta_n}]\neq0$
in $H_{\ast}(Grass_{\zeta};\Z/2\Z)$. 
Choose chambers $\si_i$ containing the $\xi_i$ in their closures.
($\si_i$ is unique if $\xi_i$ is regular.)
The choice of chambers determines cycles $\bar C_{\eta_i}$
representing the Schubert classes.
Since their homological intersection is non-trivial 
we have 
$\bar C_{\eta_1}\cap\dots\cap\bar C_{\eta_n}\neq\emptyset$.
Let $\theta$ be a point in the intersection.
All points on the Schubert {\em cell} $C_{\eta_i}$ have the same relative
position with respect to the chamber $\si_i$ 
and therefore $\tangle(\xi_i,\cdot)$ is constant along $C_{\eta_i}$,
namely equal to $\tangle(\tau_i,\eta_i)$,
cf.\ Lemma \ref{relpostitsdis}.
The semicontinuity of the Tits distance,
compare section \ref{sec:hadamard}, then implies that
$\tangle(\xi_i,\cdot)\leq\tangle(\tau_i,\eta_i)$ 
on the {\em cycle} $\bar C_{\eta_i}$.
Thus $\tangle(\xi_i,\theta)\leq\tangle(\tau_i,\eta_i)$.
It follows that 
$\sum_i m_i\cos\tangle(\tau_i,\eta_i) \leq
\sum_i m_i\cos\tangle(\xi_i,\theta)=
-slope_{\mu}(\theta)\leq0$
where $\mu$ is the measure on $\geo X$ given by the weighted
configuration $\psi$. 
Hence the inequality (\ref{stabineqlie}) holds 
whenever the corresponding Schubert classes have non-zero homological
intersection. 
\qed

\begin{rem}
\label{rem:shorterlist}
Our argument shows that 
if a semistable configuration with $\De$-weights $h$ exists 
then all inequalities hold where the homological intersection of Schubert
classes is nontrivial but not necessarily a point. 
This is an in general larger list of inequalities 
which hence has the same set of solutions in
$\De_{euc}^n$.
\end{rem}

The proof of Theorem \ref{generalschubertinequalities} works in exactly the
same way in the {\em complex} case.
The result which one obtains for complex Lie groups 
is stronger because one can work with {\em integral} cohomology and
obtains a shorter list of inequalities:

\begin{thm}
[Stability inequalities for semisimple complex groups]
\label{maincomplex}
If $G$ is complex,
then for $h\in\De_{euc}^n$ 
there exists a semistable weighted configuration with $\De$-weights $h$
if and only if 
(\ref{stabineqlie})
holds whenever the intersection 
of the integral Schubert classes 
$[C_{\eta_1}],\dots,[C_{\eta_n}]$ 
in $H_{\ast}(Grass_{\zeta};\Z/2\Z)$ 
equals $[pt]$.
\end{thm}

\subsection{The weak stability inequalities}
\label{sec:wtrineq}

In this section we consider 
a subsystem of the stability inequalities 
which corresponds to particularly simple intersections of Schubert cycles. 
We will see that it has a beautiful geometric interpretation
in terms of convex hulls.

Suppose that $G$ is a semisimple complex group.
It is well-known,
cf.\ \cite[sec.\ 2.1]{KumarLeebMillson}, 
that the Schubert cycles in the Grassmannian $Grass_{\zeta}$ 
come in pairs of mutually dual cycles, that is,
the Poincar\'e dual of a Schubert cycle is also a Schubert cycle.
The pairs of dual cycles can be nicely described 
using the parametrization of Schubert cycles 
by vertices in the Weyl group orbit $W\zeta$. 
To do so, 
let us denote by $w_0$ the element (of order two) in the Weyl group 
which maps the spherical Weyl chamber $\De_{sph}$ 
to the opposite chamber in the model apartment.
Then it has been shown in 
\cite[Lemma 2.9]{KumarLeebMillson} 
that for every vertex $\eta\in W\zeta$ holds 
\begin{equation}
\label{pdsc}
[C_{\eta}] \cdot [C_{w_0\eta}]=[pt]
\end {equation}
Hence the inequality 
(\ref{stabineqlie}) 
parametrized 
by the $n$-tuple of Schubert cycles 
$C_{\eta}$, $C_{w_0\eta}$ and 
$n-2$ times the top-dimensional cycle $C_{w_0\zeta}=Grass_{\zeta}$ 
belongs to the system of stability inequalities.
We call the subsystem consisting of all these inequalities 
for all Grassmannians $Grass_{\zeta}$ 
the {\em weak stability inequalities}.

To make them explicit,
let $h_1,\dots,h_n$ denote the $\De$-weights 
of a semistable weighted configuration on the Tits boundary 
of the symmetric space $X=G/K$,
and let $\la_{\zeta}$ denote the fundamental coweight 
which generates the edge of $\De_{euc}$ 
pointing towards the vertex $\zeta$ of $\De_{sph}$.
Then for $\eta=w\zeta$ 
the weak stability inequality 
corresponding to the intersection (\ref{pdsc}) reads:
\begin{equation*}
\langle h_1,w\la_{\zeta} \rangle + 
\langle h_2,w_0w\la_{\zeta} \rangle + 
\langle h_3,w_0\la_{\zeta} \rangle 
+\dots+
\langle h_n,w_0\la_{\zeta} \rangle 
\leq0
\end{equation*}
Using the natural isometric involution 
$h\mapsto -w_0h=:h^{\sharp}$ of $\De_{euc}$
it becomes:
\begin{equation}
\label{rewr}
\langle w^{-1}(h_1-h_2^{\sharp}),\la_{\zeta} \rangle 
\leq
\langle h_3^{\sharp}+\dots +h_n^{\sharp},\la_{\zeta} \rangle 
\end{equation}
Let 
$\De^{\ast}$
denote the (obtuse) cone 
$\{h\in{\goth a}|\langle h,\la_{\zeta}\rangle\geq0\,\forall\zeta\}$
dual to the (acute) cone $\De_{euc}$,
and let ``$\leq$'' be the order on $\mathfrak{a}$
corresponding to $\De^{\ast}$,
i.e.\ 
$h\in{\goth a}$ satisfies $h\geq 0$ if and only if $h\in\De^{\ast}$.
This order is important in representation theory.
It is usually referred to as the {\em dominance order}. 
The (sub)system of the inequalities (\ref{rewr}) 
for fixed $w$ and varying $\zeta$
amounts to the condition 
$w^{-1}(h_1-h_2^{\sharp})\in (h_3^{\sharp}+\dots+h_n^{\sharp})-\De^{\ast}$
and can hence be rewritten as the vector inequality 
\begin{equation}
\label{prewtiorder}
w^{-1}( h_1-h_2^{\sharp})\leq h_3^{\sharp}+\dots+h_n^{\sharp}
\end{equation}

\begin{thm}
[Weak stability inequalities]
Let $G$ be a semisimple complex group.
Then the $\De$-weights 
$h_1,\dots,h_n$ 
of a semistable weighted configuration
of $n$ points on the Tits boundary of the symmetric space $X=G/K$ 
satisfy for each $w \in W$ the inequality
\begin{equation}
\label{wtiord}
w h_1^{\sharp} \leq w h_2 + (h_3+\dots+h_n). 
\end{equation}
Moreover,
this system of vector inequalities 
is equivalent to the geometric condition 
\begin{equation}
\label{wtigeom}
 h_1^{\sharp} \in h_2 + convex \ hull(W\cdot(h_3+\dots+ h_n)).
\end{equation}
\end{thm}
\proof
We proved the first part already.
It follows from (\ref{prewtiorder})
by applying the order preserving involution $h\mapsto h^{\sharp}$ of
$\De_{euc}$
and renaming $w$.
(Note that 
$(wh)^{\sharp}=-w_0wh=(w_0ww_0^{-1})h^{\sharp}$.)

For the second part note that the system of inequalities (\ref{wtiord}) 
is equivalent to 
$W(h_1^{\sharp}-h_2) \subset (h_3+\dots+h_n)-\De^{\ast}$
and hence to 
\begin{equation*}
h_1^{\sharp}-h_2\in \bigcap_{w\in W} w((h_3+\dots+h_n)-\De^{\ast})=
convex \ hull (W\cdot(h_3+\dots+h_n))
\end{equation*}
A proof of the last equality 
(in the case in which $h_3 + \cdots + h_n$ is in the interior of $\De$)
can be found in \cite[pp.\ 138-140]{BorovikGelfandWhite}. 
\qed

In the case $G=GL(m,\C)$ and $n=3$
the weak stability inequalities are due to Wielandt \cite{Wielandt} 
and their geometric interpretation to Lidskii \cite{Lidskii}, 
see the discussion in the first chapter of \cite{Fulton2000}.

The simplest weak triangle inequality is inequality (\ref{wtiord})
corresponding to $w=e$,
that is to the $n$-tuples of Schubert cycles in the Grassmannians
$Grass_{\zeta}$ where one cycle is a point and the other $n-1$ are
top-dimensional:
\begin{equation}
\label{wtiordspec}
h_1^{\sharp} \leq h_2 +\dots+h_n
\end{equation}
This inequality is proven in \cite{AbelsMargulis},
compare inequality (2.28) therein.

\begin{rem}
In section \ref{sec:rank2} we will see 
for the simple complex groups of rank two 
that the weak stability inequalities are equivalent 
to the full stability inequalities,
i.e.\ all non-weak stability inequalities are redundant in these cases.
However, 
as the rank increases one would expect that most of the irredundant
inequalities are non-weak. 
Indeed, irredundant non-weak inequalities 
can already be found in rank three 
among the inequalities for the group $Sp(6,\C)$, 
see \cite[p.\ 187]{KumarLeebMillson}. 
\end{rem}
\begin{rem}
These two remarks will apply 
after section \ref{sec:pol}
where we relate the $\De$-weights of configurations on $\tits X$
with the $\De$-side lengfths of polygons in $X$. 

(i)
Notice that 
the weak stability inequalities {\em depend only on the Weyl group}
and not on further properties of the Schubert calculus,
and therefore also their solution set ${\cal W}_n(X)\subset\De_{euc}^n$.  
After proving 
Theorem \ref{polxconf:samelengths} 
and combining it with Theorem 1.3 of \cite{polbuil}
- compare our discussion in the introduction -
we will know that also the solution set 
${\cal P}_n(X)\subset\De_{euc}^n$ 
to the stability inequalities depends only on the Weyl group.
This will imply that 
${\cal P}_n(X)\subseteq{\cal W}_n(X)$,
i.e.\ 
the weak stability inequalities are a consequence 
of the stability inequalities 
not only for a semisimple complex group 
but for any noncompact semisimple real Lie group $G$.
(We are not claiming that the weak stability inequalities 
always are a subsystem of the stability inequalities
although this seems to be true also.) 

(ii)
After proving Theorem \ref{polxconf:samelengths} 
we will know that for the $\De$-side lengths $\al,\beta,\ga$ 
of an oriented geodesic triangle in $X$ 
inequality (\ref{wtiord}) 
takes the form
$w\al^{\sharp}\leq w\beta+\ga$.
Using the notation $\si(x,y)$ 
for the $\De$-distance of the oriented segment $\ol{xy}$
introduced in section \ref{sec:desidel}
we see that 
for any triple of points $x,y,z\in X$ 
and each $w\in W$ holds 
\begin{equation*}
w\cdot\si(x,z) \leq w\cdot\si(x,y) + \si(y,z).
\end{equation*}
(We use here that 
$\si(x,z)=\si(z,x)^{\sharp}$.)
The special case (\ref{wtiordspec}) turns into a nice vector valued
generalization of the ordinary triangle inequality:
\begin{equation*}
\si(x,z) \leq \si(x,y) + \si(y,z)
\end{equation*}
\end{rem}

\section{Comparing stability with Mumford stability}
\label{git}

To justify our notion of stability 
for measures on topological spherical Tits buildings,
cf.\ \cite[Definition 4.1]{polbuil} and  
Definitions \ref{def:stabmeasinf} and \ref{def:stabconfbuil},
we explain in this section
that in the example of 
weighted $n$-point configurations on complex projective space $\C\P^m$ 
our notion agrees with Mumford stability in geometric invariant theory.

Mumford introduced his notion of stability 
in order to construct good quotients 
for certain algebraic actions on projective varieties.
We start by recalling the definition 
and related concepts from geometric invariant theory, 
cf.\ \cite{Mumford,Ness}.

Consider first the case of a linear action on projective space: 
Let $V$ be a finite-dimensional complex vector space, 
let $G$ be a connected complex reductive group
and suppose that 
$\rho:G\ra SL(V)$ is a linear representation 
with finite kernel.
According to Mumford a nonzero vector $v\in V$ is called {\em unstable}
if $0\in\ol{Gv}$
and {\em semistable} otherwise. 
One obtains a notion of stability for the orbits 
of the projectivized action 
$G\acts\P(V)$. 

The Hilbert-Mumford criterion 
asserts that one can test stability on one-parameter subgroups:
$v$ is unstable if and only if there exists 
a one-parameter group $\la\subset G$
such that $0\in\ol{\la\cdot v}$,
more precisely,
if and only if there exists $\al\in{\goth p}=i{\goth k}$
such that 
$\lim_{t\to+\infty}e^{t\al}v=0$.
Let $v=\sum v_i$ be a decomposition of $v$ into eigenvectors for $\al$
and let $a_i\in\R$ be the corresponding weights,
i.e.\ $e^{t\al}v=\sum e^{a_it}v_i$.
Then 
\begin{equation}
\label{limex}
\lim_{t\to+\infty}e^{-t\mu([v],\al)} \cdot e^{t\al}v
\end{equation}
exists and is non-zero
where 
$\mu([v],\al):=\max\{a_i|v_i\neq0\}$
is a so-called {\em numerical function}. 
Hence $v$ is unstable if and only if $M([v])<0$
where $M$ denotes the derived numerical function 
$M([v]):=\inf_{0\neq\al\in i{\goth k}}\frac{\mu([v],\al)}{\|\al\|}
<0$
on $\P(V)$.

If $G\acts Y$
is an algebraic action of $G$
on an abstract projective variety $Y$
then one needs further data
in order to be able to talk about stable orbits for this action.
For instance,
it would suffice to choose 
a projective embedding $Y\subseteq\P(V)$
together with a linearization $G\ra SL(V)$ of the given action.
Then a point $[v]\in Y$, respectively, its $G$-orbit 
is called semistable 
if the vector $v$ is semistable in the above sense. 

\medskip
We restrict now to the special case
of spaces of weighted configurations on complex projective space. 
Let us consider the diagonal action of $G=SL(m+1,\C)$ 
on the projective variety
$Y\cong \times _{i=1}^n \C\P^m$
which we regard as the space of $n$-point configurations 
$\xi=(\xi_1,\dots,\xi_n)$ on $\C\P^m$. 
A choice of integral weights 
$r=(r_1,\dots,r_n)$
determines a natural projective embedding of $Y$.
Namely, 
put $W=\C^{m+1}$ and 
$W^{\otimes r}= \otimes_{i=1}^nW^{\otimes r_i}$.
The map 
$\iota:W^n \to W^{\otimes r}$ given by
\begin{equation*}
\iota(w)=
\iota (w_1,\cdots, w_n) = 
w_1^{\otimes r_1} \otimes \cdots \otimes w_n^{\otimes r_n}
\end{equation*}
induces the Segre embedding
$Y\embed \P(W^{\otimes r})$. 
The natural $G$-action on $W^{\otimes r}$ 
linearizes the given action $G\acts Y$. 
The configuration 
$([w_1],\dots,[w_n])\in Y$ 
is Mumford semistable 
(with respect to the chosen embedding and linearization) 
if and only if the $G$-orbit
\[
g\mapsto 
g (w_1^{\otimes r_1} \otimes \cdots \otimes w_n^{\otimes r_n})
=
(g w_1)^{\otimes r_1} \otimes \cdots \otimes (g w_n)^{\otimes r_n}
\]
does not accumulate at $0$,
that is, 
if and only if the orbital distance function 
\[ \psi_w(g) := 
\|g (w_1^{\otimes r_1} \otimes \cdots \otimes w_n^{\otimes r_n}) \|
=
\|gw_1\|^{r_1}\cdots \|gw_n\|^{r_n} \]
is bounded away from zero.
Here lengths are measured with respect to a fixed Hermitian form on $W$
and the induced Hermitian form on $W^{\otimes r}$. 

Let us compare this with our notion of stability.
We may view 
$\xi=([w_1],\dots,[w_n])$
as a weighted configuration on the ideal boundary 
of the symmetric space $X=G/K$,
$K=SU(m+1)$, 
since $\C\P^m$ canonically identifies 
with a maximally singular $G$-orbit on $\geo X$.
Moreover, 
we may choose the norm $\|\cdot\|$ on $W$ to be $K$-invariant.
The connection between the orbital distance function $\psi_w$ 
and weighted Busemann functions on $X$ is based on the fact 
- cf.\ Example \ref{ex:busesl} - 
that after suitable normalization of the metric on $X$ 
we have 
$\log\|g^{-1}w_j\| = b_{[w_j]}(gK)$ 
modulo additive constants. 
Therefore 
\begin{equation*}
\log\psi_w(g^{-1}) + const 
= \sum_{j=1}^n r_j\cdot b_{[w_j]}(gK) 
= b_{\mu}(gK)  
\end{equation*}
where 
$\mu=\sum_{j=1}^n r_j\cdot \de_{[w_j]}$
is the measure associated to the weighted configuration $\xi$.
Thus the configuration $\xi$ is semistable in Mumford's sense if and only if  
the weighted Busemann function $b_{\mu}$ is bounded below.
If $b_{\mu}$ is bounded below
then $slope_{\mu}\geq0$ on $\geo X$,
that is, $\xi$ is semistable in our sense,
cf.\ Definition \ref{def:stabmeasinf}.

Both notions of (semi)stability are in fact equivalent
in this case. 
It is not hard to show this directly.
It also follows from the Hilbert-Mumford criterion 
together with the observation that 
the numerical function $\mu(\xi,\cdot)$
is essentially the slope function $slope_{\mu}$.
Namely,
the existence of the limit (\ref{limex}) implies that 
\begin{equation*}
b_{\mu}(e^{t\al}K)
=\log\psi_w(e^{-t\al})+const
=\mu(\xi,-\al)\cdot t+O(1)
\end{equation*}
where $\mu(\xi,\al)$ is the maximal weight $a_i$ 
for the action of $\al\in i\cdot {\goth s \goth u}(m+1)$ on $W^{\otimes r}$
such that the corresponding component of $w^{\otimes r}$ is non-zero.
The lines $t\mapsto e^{t\al}K$ 
are the geodesics in $X$ through the base point $o$ fixed by $K$.

\section{Relating polygons and configurations}
\label{sec:pol}

\subsection{The $\De$-side lengths of oriented polygons in $X$ and ${\goth p}$}
\label{sec:desidel}

The equivalence classes of oriented geodesic segments 
in the symmetric space $X=G/K$ modulo the natural $G$-action by isometries 
are parametrized by the Euclidean Weyl chamber $\De_{euc}$ 
associated to $X$,
\[ G\backslash X\times X \cong \De_{euc} .\]
The vector $\si(x,y)\in \De_{euc}$ corresponding to an oriented segment $\ol{xy}$ 
can hence be thought of as a vector-valued length.
We call it the {\em $\De$-length} of the oriented segment. 

We think of $X^{\Z/n\Z}$, $n\geq3$, 
as the space of oriented closed $n$-gons in $X$.
An $n$-tuple $(x_1,\dots,x_n)$ 
is interpreted as the polygon with vertices
$x_1,\dots,x_n$ and the $i$-th edge $\ol{x_{i-1}x_i}$ is denoted by $e_i$. 
(Recall that any two points in $X$ are connected by a unique geodesic segment.)
In the sequel,
all polygons will be assumed to be {\em oriented}. 

We denote by 
\begin{equation*}
\label{sidelengthmap}
\si:X^{\Z/n\Z}\lra \De_{euc}^n, \quad
(x_1,\dots,x_n)\mapsto(\si(x_0,x_1),\dots,\si(x_{n-1},x_n))
\end{equation*}
the side length map.
We are interested in its image 
\begin{equation*}
{\cal P}_n(X) \subset \De_{euc}^n .
\end{equation*}

We are also interested in the 
analogous problems for the infinitesimal symmetric space 
$T_oX\cong{\goth p}$ associated to $G$ 
- compare the terminology and notation from section \ref{sec:infsymm} - 
namely to study the image 
\begin{equation*}
{\cal P}_n({\goth p}) \subset \De_{euc}^n 
\end{equation*}
of the natural $\De$-side length map 
$\si':{\goth p}^{\Z/n\Z} \lra \De_{euc}^n$. 

To solve both problems,
we will now translate them into a question about weighted configurations at infinity 
which has been treated in section \ref{sec:conf}. 
Namely, we will prove that the $\De$-side length spaces 
${\cal P}_n(X)$ and ${\cal P}_n({\goth p})$
coincide with the space of possible weights of semistable configurations on $\tits X$.
The relation between polygons in $X$, respectively, in ${\goth p}$
and weighted configurations on $\tits X$
is established by a Gauss map type construction 
as in \cite[sec.\ 4.2]{polbuil},
respectively, by radial projection.

\subsection{Relating polygons in ${\goth p}$ and configurations on $\tits X$}
\label{sec:polyconf}

An $n$-tuple
$e=(e_1,\dots,e_n)\in{\goth p}^n$ 
can be interpreted 
as an {\em open $n$-gon} in ${\goth p}$ 
with vertices 
$v_i=\sum_{j\leq i}e_j$, $0\leq i\leq n$ 
and edges $e_i$. 
The $n$-gon closes up 
if and only if the {\em closing condition}
\begin{equation}
\label{closecond}
e_1+\dots+e_n=0
\end{equation}
holds.
(The distinction between open and closed $n$-gons 
will be restricted to this section.
Elsewhere in this paper $n$-gons are supposed to be closed.)

An open $n$-gon $e$ corresponds 
to a weighted $n$-point configuration $\psi$ on $\tits X$
by assigning to each non-zero edge $e_i$ the mass $m_i:=\|e_i\|$ 
located at the ideal point $\xi_i$ corresponding to $e_i$
under the radial projection ${\goth p}-\{0\}\ra\geo X$. 
In the case $e_i=0$ we set $m_i=0$ and choose $\xi_i$ arbitrarily.
Note that 
the $\De$-weights of $\psi$ equal the $\De$-side lengths of $e$.
This is what makes this correspondence between 
polygons and weighted configurations 
useful for us. 

For a unit vector $v\in {\goth p}$ and the corresponding ideal point $\eta\in\geo X$ 
holds $\nabla b_{\eta}(o)=-v$.
Hence 
\begin{equation*}
e_1+\dots+e_n=-\nabla b_{\mu}(o) .
\end{equation*}
The closing condition (\ref{closecond})
is therefore equivalent to $o$ being a minimum of $b_{\mu}$.
(Since $b_{\mu}$ is a convex function,
its critical points are global minima.) 
This implies: 

\begin{lem}
The weighted configurations on $\tits X$ 
corresponding to closed polygons in ${\goth p}$
are nice semistable. 
\end{lem}

Conversely,
if $\psi$ is nice semistable with $\De$-weights $h$,
we may use the natural $G$-action on weighted configurations 
to move a critical point of $b_{\mu}$ 
- which exists by Lemma \ref{nss} - 
to the base point $o$.
Thus:

\begin{lem}
\label{gnssfrompol}
Let $\psi$ be a nice semistable weighted configuration on $\tits X$.
Then its $G$-orbit contains a configuration which corresponds 
to a closed polygon. 
\end{lem}

Given a semistable configuration,
Lemma \ref{nssinclosure} tells 
that the closure of its $G$-orbit contains a nice semistable configuration 
which then has the same $\De$-weights.
Hence the $\De$-weights of semistable configurations 
occur also for nice semistable configurations
and we conclude: 

\begin{thm}
\label{equivalenceofmoduli}
For $h\in\De_{euc}^n$ 
there exist closed $n$-gons in ${\goth p}$ with $\De$-side lengths $h$ 
if and only if there exist 
semistable weighted configurations on $\tits X$ with $\De$-weights $h$.
\end{thm}

\medskip
Let us briefly specialize to the case
when $G$ is a {\em complex} semisimple Lie group. 
We then have ${\goth p}=i{\goth k}$ 
and may identify ${\goth k}$ and ${\goth p}$ as $K$-modules.
$K$ acts on ${\goth k}$ by the adjoint action. 

Given $h=(h_1,\dots,h_n)\in\De_{euc}^n$
we let ${\cal O}_i$ denote the $K$-orbit in ${\goth k}$ corresponding to $h_i$.
The product space 
${\cal O}_1\times\dots\times{\cal O}_n$ 
is naturally identified with the space of open $n$-gons in ${\goth p}$ 
with fixed $\De$-side lengths $h$.
All $K$-orbits ${\cal O}$ in ${\goth k}$ carry natural invariant symplectic structures. 
It is a standard fact that the momentum maps 
for the actions $K\acts{\cal O}$ are given by the embeddings 
${\cal O}\embed{\goth k}$
where one identifies ${\goth k}^*\cong{\goth k}$
via the Killing form. 
Hence:
\begin{lem}
\label{momentdiagaction}
The diagonal action 
$K\acts {\cal O}_1\times\dots\times{\cal O}_n$
is Hamiltonian with momentum map 
${\cal O}_1\times\dots\times{\cal O}_n\ra{\goth k}\cong{\goth p}$
given by 
\[ m(e)=\sum_{i=1}^n e_i. \]
\end{lem}
We see that 
the closing condition (\ref{closecond}) amounts in this situation
to the momentum zero condition from symplectic geometry,
cf.\ \cite{Kirwan}.

\subsection{Relating polygons in $X$ and configurations on $\tits X$}
\label{sec:polyxconf}

We prove results analogous to those in section \ref{sec:polyconf}
for polygons in the symmetric space $X=G/K$.
They are more difficult to obtain 
because the closing condition is nonabelian
and not directly related to the differential 
of the weighted Busemann function.  
We consider only closed polygons. 

Let $P$ be a $n$-gon in $X$, 
i.e.\ a map $P:\Z/n\Z\ra X,i\mapsto x_i$.
Its side lengths $m_i=d(x_{i-1},x_i)$
determine a measure $\nu$ on $\Z/n\Z$ by putting $\nu(i)=m_i$.
and $P$ gives rise to a {\em Gauss map} 
\begin{equation*}
\psi:\Z/n\Z \lra \tits X
\end{equation*}
by assigning to $i$ an ideal point $\xi_i\in\tits X$
so that the ray $\ol{x_{i-1}\xi_i}$ passes through $x_i$;
the point $\xi_i$ is unique unless $x_{i-1}=x_i$.
This construction, in the case of hyperbolic plane, 
already appears in the letter of Gauss to W.\ Bolyai \cite{Gauss}.
Taking into account the measure $\nu$,
we view $\psi$ as a weighted configuration on $\tits X$.
The $\De$-weights of $\psi$ equal the $\De$-side lengths of the polygon $P$. 

The next result shows that 
again the configurations arising from polygons are characterized 
by a stability property,
namely they are nice semistable.
That they are semistable
is proven in Lemma 4.3 of \cite{polbuil} in a more general situation
(and this is actually all we need in the proof of our main results,
cf.\ Theorem \ref{polxconf:samelengths}). 

\begin{lem}
[Nice semistability of Gauss map]
\label{gaussss}
The weighted configurations on $\tits X$ 
arising as Gauss maps of closed polygons in $X$
are nice semistable. 
\end{lem}
\proof
For the convenience of the reader 
we first reproduce the proof of 
Lemma 4.3 of \cite{polbuil} 
to show semistability. 
Let $\ga_i:[0,m_i]\ra X$ be unit speed parametrizations 
of the sides $\ol{x_{i-1}x_i}$ of the polygon $P$
and let $\eta\in\tits X$ be an ideal point. 
The derivative of the Busemann function $b_{\eta}$ along $\ga_i$
is given by 
\begin{equation}
\label{derest}
\frac{d}{dt}
(b_{\eta}\circ\ga_i)(t)
= -\cos\angle_{\ga_i(t)}(\xi_i,\eta)
\leq -\cos\tangle(\xi_i,\eta) ,
\end{equation}
compare formula (\ref{diffbuse}) in section \ref{sec:symmsp}.
Integrating along $\ga_i$ we obtain
\begin{equation*}
b_{\eta}(x_i)-b_{\eta}(x_{i-1}) \leq
-m_i\cdot\cos\tangle(\xi_i,\eta) 
\end{equation*}
and summation over all sides yields
\[ 0 \leq -\sum_{i\in\Z/n\Z} m_i\cdot\cos\tangle(\xi_i,\eta) = slope_{\mu}(\eta) .\]
confirming the semistability of the measure $\mu=\psi_*\nu$
and the configuration $\psi$.

Regarding {\em nice} semistability,
suppose that $\psi$ is not stable,
i.e.\ that $S:=\{slope_{\mu}=0\}$ is non-empty. 
For an ideal point $\eta\in S$
we have equality in (\ref{derest}),
that is, 
$b_{\eta}$ is linear along every segment $\ga_i$.
Denote by $l_i$ the line passing through $x_i$ and asymptotic to $\eta$. 
Lemma \ref{atomicbusemannflinearalonggeo} implies for $i=1,\dots,n$
that the two lines $l_{i-1}$ and $l_i$ are parallel.
It follows that the polygon $P$ is contained in the parallel set 
$P(l_1)=\dots=P(l_n)$.
Moreover, $\mu$ is supported on its ideal boundary, 
We denote by $\hat\eta$ the other ideal endpoint of the lines $l_i$. 
Since $\mu$ is supported on $\geo P(l_i)$, 
the Busemann function $b_{\mu}$ is linear on all lines (parallel to) $l_i$ 
and hence 
$slope_{\mu}(\hat\eta)=-slope_{\mu}(\eta)=0$,
i.e.\ $\hat\eta\in S$.

Recall that $S$ is a convex subcomplex of $\tits X$ 
because $\mu$ is semistable.  
We may choose $\eta$ maximally regular in $S$,
that is, as an interior point of a top-dimensional simplex $\si\subset S$.
Then $S$ contains the convex hull $s$
of $\si$ and $\hat\eta$
which is a top-dimensional unit sphere in $S$, $\dim(s)=\dim(\si)=\dim(S)$.
Thus $\mu$ is nice semistable
according to Definition \ref{def:nss}.
\qed

\medskip
We are now interested in finding polygons with prescribed Gauss map.
Such polygons will correspond to the fixed points 
of a certain weakly contracting self map of $X$,
compare section 4.3 in \cite{polbuil}.

For $\xi\in \tits X$ and $t\geq0$,
we define the map 
$\phi_{\xi,t}: X\to X$
by sending $x$ to the point at distance $t$ from $x$
on the geodesic ray $\ol{x\xi}$.
Since $X$ is nonpositively curved,
the function 
$\delta: t\mapsto d(\phi_{\xi,t}(x), \phi_{\xi,t}(y))$
is convex.
It is also bounded because the rays $\ol{x\xi}$ and $\ol{y\xi}$ are asymptotic, 
and hence it is monotonically non-increasing in $t$.
This means that the maps $\phi_{\xi,t}$ are weakly contracting, 
i.e.\ they are $1$-Lipschitz.
For the weighted configuration $\psi$ we define the weak contraction 
\begin{equation*}
\Phi_{\psi}:X\lra X
\end{equation*}
as the composition $\phi_{\xi_n,m_n}\circ\dots\circ\phi_{\xi_1,m_1}$.
The fixed points of $\Phi_{\psi}$
are the $n$-th vertices of polygons $P=x_1\ldots x_n$
with Gauss map $\psi$.
The next result is the counterpart of Lemma \ref{gnssfrompol}.

\begin{prop}
\label{prop:stabletofixed}
If the weighted configuration $\psi$ on $\tits X$ is nice semistable 
then the weak contraction $\Phi_\psi:X\to X$ has a fixed point.
\end{prop}

\proof
We will use the following auxiliary result
which extends Cartan's fixed point theorem 
for isometric actions on Hadamard spaces with bounded orbits. 

\begin{lem}[{\cite[Lemma 4.6]{polbuil}}]
\label{bounded}
Let $Y$ be a Hadamard space and $\Phi: Y\to Y$ a 1-Lipschitz self map.
If the forward orbits $(\Phi^n y)_{n\ge 0}$ are bounded
then $\Phi$ has a fixed point in $Y$.
\end{lem}

It therefore suffices to show
that the dynamical system $\Phi_{\psi}:X\ra X$ has a bounded forward orbit
$(\Phi_{\psi}^n (p))_{n\ge 0}$.
Suppose that this is false.

{\em Step 1.} Our assumption that $\Phi_{\psi}$ 
does not have a bounded forward orbit
implies that $\Phi_{\psi}$ does not map any bounded subset of $X$ into itself.
Pick a base point $o\in X$. Since no metric ball centered at $o$ is
mapped into itself there is a sequence of points $x_n$ with $d(x_n,o)\ra\infty$
which is ``pulled away'' from $o$ in the sense that
\[  d(\Phi_{\psi}(x_n),o) > d(x_n,o) . \]
Since $\Phi_{\psi}$ is 1-Lipschitz
we have in fact 
$d(\Phi_{\psi}(x),o) > d(x,o)$
for all points $x$ on all segments $\ol{ox_n}$.
Since $X$ is locally compact,
after passing to a subsequence, 
the segments $\ol{ox_n}$ 
Hausdorff converge to a geodesic ray 
$\rho([0,+\infty))=\ol{o \xi}$.
This ray is ``pulled away'' from $o$ in the sense that for each $t\geq0$ holds 
\begin{equation}
\label{pulloff}
d(\Phi_{\psi}(\rho(t)),o) \ge d(\rho(t),o).
\end{equation}

{\em Step 2.}
For any unit speed geodesic ray $\rho_0:[0,+\infty)\to X$ 
we claim that 
\begin{equation}
\label{pullray}
\lim_{t\to+\infty} d(\Phi_{\psi}(\rho_0(t)),o)-d(\rho_0(t),o)
= -slope_{\mu}(\eta)
\end{equation}
where $\eta$ is the ideal endpoint of $\rho_0$
and $\mu=\psi_*\nu$ the measure associated 
to the weighted configuration $\psi$. 
To verify this
we first look for a ray $\rho_1$ such that 
$d(\phi_{\xi_1,m_1}(\rho_0(t)),\rho_1(t))$ $\to 0$
for $t\to+\infty$.
(Note that 
$\phi_{\xi_1,m_1}\circ\rho_0$ is in general no geodesic ray.) 
There exists a geodesic line $l_1$ asymptotic to $\eta$ 
such that $\xi_1\in\geo P(l_1)$. 
Inside the parallel set $P(l_1)$ 
there is a unique ray $\hat\rho_0$ strongly asymptotic to $\rho_0$.
The weak contraction $\phi_{\xi_1,m_1}$ 
maps lines parallel to $l_1$ again to such lines 
and 
$\rho_1:=\phi_{\xi_1,m_1}\circ\hat\rho_0$ 
is a geodesic ray with the desired property.
Since 
\begin{equation}
\label{klbslbr}
b_{\eta}\circ\phi_{\xi_1,m_1}-b_{\eta}\equiv-m_1\cdot\cos\tangle(\eta,\xi_1)
\qquad\hbox{on $P(l_1)$}
\end{equation}
we have 
$b_{\eta}\circ\rho_1(t)-b_{\eta}\circ\rho_0(t)
\to-m_1\cdot\cos\tangle(\eta,\xi_1)$
for $t\to+\infty$.

Proceeding by induction,
we find rays $\rho_1,\dots,\rho_n$ asymptotic to $\eta$
such that for $i=1,\dots,n$ holds 
\begin{equation*}
d(\phi_{\xi_i,m_i}(\rho_{i-1}(t)),\rho_i(t)) \to 0
\end{equation*}
and 
\begin{equation*}
b_{\eta}\circ\rho_i(t)-b_{\eta}\circ\rho_{i-1}(t)
\to
-m_i\cdot\cos\tangle(\eta,\xi_i) 
\end{equation*}
as $t\to+\infty$.
It follows using the weak contraction property that 
\begin{equation*}
d(\Phi_{\psi}(\rho_0(t)),\rho_n(t)) \to 0
\end{equation*}
and 
\begin{equation*}
b_{\eta}\circ\rho_n(t)-b_{\eta}\circ\rho_0(t)
\to slope_{\mu}(\eta)
\end{equation*}
as $t\to+\infty$.
Note that with the normalization $b_{\mu}(o)=0$ 
any ray $\tilde\rho$ asymptotic to $\eta$ satisfies 
$d(o,\tilde\rho(t))+b_{\mu}(\tilde\rho(t))\to0$
as $t\to+\infty$.
As a consequence we obtain (\ref{pullray}).

{\em Step 3.}
Suppose first that the configuration $\psi$ is {\em stable}. 
Choosing $\rho_0=\rho$ 
we obtain in view of $slope_{\mu}(\xi)>0$
a contradiction between (\ref{pullray}) and (\ref{pulloff}).

We are left with the case 
that the configuration $\psi$ is nice semistable but {\em not stable}.
According to Definition \ref{def:nss},
the $d$-dimensional convex subcomplex $\{slope_{\mu}=0\}$ of $\tits X$
contains a unit $d$-sphere $s$.
We consider a $(d+1)$-flat $f\subset X$
such that $\geo f=s$
and inside $f$ a maximally regular geodesic line $l$.
Then $P(f)=P(l)$.
Since $b_{\mu}$ is constant on $f$,
Lemma \ref{weightedbusemannflinearalonggeo} implies 
that the measure $\mu$ is supported on $\geo P(f)$. 

For any geodesic $l'\subset f$
we therefore have that $\mu$ is supported on $\geo P(l')\supseteq\geo P(f)$.
Denoting the two ideal endpoints of $l'$ in $s$ by $\eta'_{\pm}$ 
we obtain as in (\ref{klbslbr}) that 
$b_{\eta'_{\pm}}\circ\Phi_{\psi}-b_{\eta'_{\pm}}
\equiv slope_{\mu}(\eta'_{\pm})=0$ on $P(l')$,
that is,
$\Phi_{\psi}$ preserves each cross section $\{pt\}\times CS(l')$
of $P(l')\cong l'\times CS(l')$.
Since this holds for all geodesics $l'$ in $f$
we conclude that 
$\Phi_{\psi}$ preserves each cross section $\{pt\}\times CS(f)$
of $P(f)\cong f\times CS(f)$.

Let $Z$ be one of the cross sections.
As in step 1,
there exists a ray $\rho$ in $Z$
satisfying (\ref{pulloff}).
Note that $slope_{\mu}>0$ on $\geo Z$. 
Otherwise $\{slope_{\mu}=0\}$ 
would contain a hemisphere of dimension $dim(s)+1$
which is impossible because $dim(\{slope_{\mu}=0\})=dim(s)$,
compare the proof of Lemma \ref{nss}.
As in the case when $\mu$ is stable,
see the beginning of this step,
we obtain a contradiction. 
This concludes the proof of the Proposition.
\qed

\begin{rem}
If the configuration $\psi$ is stable
then $\Phi_{\psi}$ has a {\em unique} fixed point.
Namely, if $x$ and $x'$ are different fixed points on a line $l$
then $\Phi_{\psi}$ restricts on $l$ to an isometry.
It follows that $\mu$ is supported on $\geo P(l)$
and hence not stable.
\end{rem}

Analogous to Theorem \ref{equivalenceofmoduli} we obtain:

\begin{thm}
\label{polxconf:samelengths}
For $h\in\De_{euc}^n$ 
there exist closed $n$-gons in $X$ with $\De$-side lengths $h$ 
if and only if there exist 
semistable weighted configurations on $\tits X$ with $\De$-weights $h$.
\end{thm}
\proof
This is a direct consequence of Lemma \ref{gaussss},
Proposition \ref{prop:stabletofixed}
and the fact already used above 
that the existence of semistable configurations on $\tits X$
with $\De$-weights $h$
implies the existence of nice semistable ones,
cf.\ Lemma \ref{nssinclosure}. 
\qed

\medskip
Combining Theorems \ref{equivalenceofmoduli} and \ref{polxconf:samelengths}
we obtain 
Theorem \ref{thm:thompson} stated in the introduction.
It generalizes the Thompson Conjecture \cite{Thompson}
which was formulated for the case of $G=GL(n,\C)$. 
Special cases were obtained in \cite{Klyachko2} and 
\cite{AMW}. 
Another proof in the general case 
has recently been given in \cite{EvensLu}.

\section{The stable measures 
on the Grassmannians of the classical groups}
\label{sec:expl}

In section \ref{sec:defstab} 
we defined stability 
for measures of finite total mass 
on the ideal boundary of symmetric spaces of noncompact type.
In this section 
we will make the stability condition explicit 
in the case of the classical groups. 
We will restrict ourselves to measures supported on a maximally singular
orbit,
that is, 
to measures supported 
on the (generalized) Grassmannians associated to the classical groups.

\subsection{The special linear groups} 

Let $G=SL(n,\C)$ and 
let $X=SL(n,\C)/SU(n)$ be the associated symmetric space. 
Recall that the Tits boundary $\tits X$, as a spherical building,
is combinatorially equivalent to the complex of flags 
of proper non-zero linear subspaces of $\C^n$. 
The vertices of $\tits X$ correspond to the linear subspaces 
and the simplices to the partial flags of such,
compare \cite[p.\ 120]{Brown}.

Let $\mu$ be a measure with finite total mass 
supported on an orbit of vertices $G\eta$ in $\tits X$.
Such an orbit 
is identified with the Grassmannian $G_q(\C^n)$ of $q$-planes 
for some $q$, $1\leq q\leq n-1$. 
(This identification is a homeomorphism with respect to the cone topology,
cf.\ section \ref{sec:symmsp}.)
The main issue in determining $slope_{\mu}$
and evaluating the system of inequalities 
(\ref{ssineq}) given in section \ref{sec:defstab} 
is to compute the distances between vertices in $\tits X$.
We denote by $[U]$ the vertex in $\tits X$ corresponding to the subspace $U$,
and for non-trivial linear subspaces $U,V\subset\C^n$
we introduce the auxiliary function
\[ \dim_U(V):=\frac{\dim(U\cap V)}{\dim(U)} \qquad.\]

\begin{lem}
The distance between the vertices of $\tits X$ 
corresponding to proper non-zero linear subspaces $U,V\subset\C^n$
is given by 
\begin{equation}
\label{titsdistofpointsingrass}
\cos\tangle([U],[V])=
C\cdot (\dim_U(V)-\frac{1}{n}\dim(V))
\end{equation}
where $C=C(\dim(U),\dim(V),n)$
is a {\em positive} constant.  
\end{lem}
\proof
Let $p=\dim(U)$, $q=\dim(V)$, $s=\dim(U\cap V)$ 
and choose a basis $e_1\dots e_n$ of $\C^n$ 
so that $e_1\dots e_p$ is a basis of $U$ and 
$e_{p-s+1}\dots e_{p+q-s}$ is a basis of $V$. 
The splitting 
$\C^n=\langle e_1\rangle \oplus\dots\oplus\langle e_n\rangle$ 
determines a maximal flat $F$ in $X$ (of dimension $n-1$) 
in the sense that 
the vertices of the apartment $\tits F\subset\tits X$ 
correspond to the non-trivial linear subspaces of $\C^n$ 
spanned by some of the $e_i$. 
Let us denote by $\hat e_i$ the unit vector field in $F$ pointing 
towards $[\langle e_i\rangle]\in\tits F$. 
Then, since $\hat e_1+\dots+\hat e_n =0$, we find by symmetry that 
$\hat e_i\cdot \hat e_j=-\frac{1}{n-1}$ for $i\neq j$. 
Note that the vector field $\hat e_{i_1}+\dots+\hat e_{i_k}$ 
points towards the vertex 
$[\langle e_{i_1}\dots e_{i_k}\rangle]\in\tits F$. 
It follows that the Tits distance between the vertices $[U]$ and $[V]$ 
is given by the angle between the vector fields 
$\hat e_1+\dots+\hat e_p$ and 
$\hat e_{p-s+1}+\dots+\hat e_{p+q-s}$ 
whose cosine equals 
\[ 
 \sqrt{\frac{n-1}{p(n-p)}} \sqrt{\frac{n-1}{q(n-q)}}
 \  (pq\frac{-1}{n-1}+s\frac{n}{n-1}) 
= const(p,q,n)\cdot(\frac{s}{p}-\frac{q}{n})
\]
whence (\ref{titsdistofpointsingrass}). 
\qed

\medskip\noindent
The slope function $slope_{\mu}$ then takes 
the following form, cf.\ (\ref{slopeform}): 
\[ slope_{\mu}([U]) = 
-\int_{G_q(\C^n)} \cos\tangle([U],[V]) \;d\mu([V]) 
.\]
Hence $slope_{\mu}([U])\geq0$ if and only if 
\begin{equation}
\label{SStabkritGrassmannianSLn}
\int_{G_q(\C^n)} \dim_U(V) \;d\mu([V]) \le \frac{q}{n} \|\mu\|. 
\end{equation}
We conclude using Corollary \ref{vertexslopes}:

\begin{prop}
\label{stabforsl}
A measure $\mu$ with finite total mass 
on the Grassmannian $G_q(\C^n)$ 
is semistable if and only if 
(\ref{SStabkritGrassmannianSLn}) 
holds for all proper non-zero linear subspaces 
$U\subset\C^n$. 
It is stable if and only if all inequalities hold strictly.
\end{prop}
\begin{ex}
Let $\mu$ be a measure with finite total mass on 
complex projective space $\C\P^m=G_1(\C^{m+1})$.
Then $\mu$ is semistable if and only if  
each $d$-dimensional projective subspace carries at most 
$\frac{d+1}{m+1}$ times the total mass. 
For instance,
a measure on the complex projective line is semistable 
if and only if each atom carries at most half of the total mass. 
\end{ex}

\subsection{The orthogonal and symplectic groups}

Let us now consider 
the other families 
of classical groups $G=SO(n,\C)$ and $G=Sp(2n,\C)$. 
To streamline the notation we note that in either case 
$G$ is the group preserving 
a non-degenerate bilinear form $b$ 
on a finite-dimensional complex vector space $W$
and a complex volume form.
We denote by $Y=G/K$ the symmetric space associated to $G$. 
The spherical Tits building $\tits Y$ is combinatorially equivalent to the 
flag complex of non-zero $b$-isotropic subspaces of $W$,
cf.\ \cite[pp.\ 123-126]{Brown}. 
(In the case of $SO(2m,\C)$ one may prefer to consider the natural 
{\em thick} building structure on $\tits Y$ combinatorially equivalent 
to the flag complex for the ``oriflamme geometry'', 
but for our considerations this does not make a difference.)

As before in the case of the special linear group 
the main task is to compute the distances between vertices of $\tits Y$.
For an isotropic subspace $U\subset W$ 
we let $(U)$ denote the corresponding vertex in $\tits Y$. 

\begin{lem}
The distance between the vertices of $\tits Y$
corresponding to non-zero 
$b$-isotropic subspaces $U,V\subset W$ 
is given by 
\begin{equation*} 
\cos\tangle((U),(V)) = C \cdot (dim_U(V) + dim_U(V^{\perp}) -1 )
\end{equation*}
where $C=C(dim(U),dim(V),dim(W))$ is positive constant. 
\end{lem}
\proof
We will use the inclusion of $G$ into the appropriate 
special linear group. 
It induces an isometric embedding
$\iota:\tits Y \embed \tits X$
of Tits boundaries,
\[ \cos\tangle((U),(V)) = \cos\tangle(\iota((U)),\iota((V))) .\]
Recall that the 
form $b$ induces an isometric automorphism (polarity) of $\tits X$ 
which we will also denote $b$.
It satisfies $b([U]) = [U^{\perp}]$. 
The point $\iota((U))$ is the midpoint of the edge in $\tits X$ 
joining the vertices $[U]$ and $[U^{\perp}]$.
(These are interchanged by $b$, and they are connected by an edge
because $U$ is isotropic.)  

There exists an apartment $a\subset\tits X$ 
containing the edge joining $[U]$ and $[U^{\perp}]$ 
and the edge joining $[V]$ and $[V^{\perp}]$. 
(We allow the possibility that $U=U^{\perp}$ or $V=V^{\perp}$.) 
Note that $\tangle([U],[U^{\perp}])$ and $\tangle([V],[V^{\perp}])$ 
depend only on $dim(W)$ and $dim(U)$ resp.\ $dim(V)$,
cf.\ (\ref{titsdistofpointsingrass}).
Hence 
\begin{eqnarray*} 
\cos\tangle(\iota((U)),\iota((V))) = 
C\cdot (\cos\tangle([U],[V])+\cos\tangle([U],[V^{\perp}])+ \\
\cos\tangle([U^{\perp}],[V])+\cos\tangle([U^{\perp}],[V^{\perp}])) 
\end{eqnarray*}
with a constant $C=C(dim(U),dim(V),dim(W))>0$.
Since $b$ is isometric,
we have 
$\cos\tangle([U^{\perp}],[V])=\cos\tangle([U],[V^{\perp}])$
and 
$\cos\tangle([U^{\perp}],[V^{\perp}])=\cos\tangle([U],[V])$. 
So 
\[ \cos\tangle(\iota((U)),\iota((V))) = 
2C\cdot (\cos\tangle([U],[V])+\cos\tangle([U],[V^{\perp}])) \]
We now can appeal to our computation of Tits distances for the special linear group 
and deduce the assertion from (\ref{titsdistofpointsingrass}).
\qed

\medskip
The Grassmannian $G^o_q(W,b)$ of $b$-isotropic $q$-planes in $W$ 
embeds as a maximally singular $G$-orbit in $\tits Y$. 
Let $\mu$ be a measure on $G^o_q(W,b)$ with finite total mass. 
The slope function $slope_{\mu}$ takes on vertices of $\tits Y$ the form 
\begin{equation*}
slope_{\mu}((U)) = 
  - C\cdot \int_{G^o_q(W,b)} (dim_U(V) + dim_U(V^{\perp}) -1 ) \, d\mu((V)) 
\end{equation*}
with $C=C(dim(U),q,dim(W))>0$. 
Using Corollary \ref{vertexslopes}, 
we obtain:
\begin{prop}
\label{stabforisotropic}
Let $b$ be a non-degenerate symmetric or alternating form on a 
finite-dimensional complex vector space $W$ 
and let $G^o_q(W,b)$ be the Grassmannian of $b$-isotropic
$q$-planes. 
A measure $\mu$ with finite total mass 
supported on $G^o_q(W,b)$ 
is semistable if and only if 
for every non-zero $b$-isotropic subspace $U\subset W$ holds 
\begin{equation*}
\int_{G^o_q(W,b)} (\dim_U(V) + dim_U(V^{\perp})) \, d\mu([V]) \leq  \|\mu\| .
\end{equation*}
It is stable if and only if all inequalities hold strictly.
\end{prop}

\section{The $\De$-side length polyhedra for the rank two root systems}
\label{sec:rank2}

In this section,
we make the stability inequalities 
given in Theorem \ref{generalschubertinequalitiesintro}
explicit for the simple complex Lie groups of rank two. 
Since the $\De$-side length polyhedron
depends only on the spherical Coxeter complex, respectively,
on the root system ${\cal R}$,
cf.\ Theorem \ref{polddepsphcox}, 
we will also denote it by ${\cal P}_n({\cal R})$.
For the root system ${\cal R}$ 
there are three possible cases,
$A_2$, $B_2=C_2$ and $G_2$,
and the corresponding simple complex Lie algebras ${\goth g}$ are 
$sl(3,\C), so(5,\C)$ and $g_2$.

For $n=3$,
i.e.\ in the case of triangles, 
we shall see that 
the systems of inequalities for $B_2=C_2$ 
and $G_2$ are redundant.  
Using the computer program Porta
we will describe the irredundant subsystems. 
It turns out that in all rank two cases 
the redundant inequalities are precisely 
the non-weak triangle inequalities 
in the sense of section \ref{sec:wtrineq}.
We will also give generators for the polyhedral cones 
${\cal P}_3({\cal R})$ 
(again using Porta).

\subsection{Setting up the notation}

Let $G$ be a simple complex Lie group with finite center 
and Lie algebra isomorphic to ${\goth g}$. 

For the root systems 
$A_2$, $B_2=C_2$ and $G_2$
the Weyl group $W$ is a dihedral group  of order $6$, $8$ and $12$ 
respectively. 
In each case there are two (generalized) Grassmannians $G/P_i$ 
associated to $G$.
They correspond to the two vertices $\zeta_1$ and $\zeta_2$ 
of the spherical Weyl chamber (arc) $\De_{sph}$. 
After identifying the spherical model Weyl chamber $\De_{sph}$
with a chamber in the Tits boundary $\tits X$ 
of the associated symmetric space $X=G/K$ 
(which is the spherical Tits building attached to $G$)
we may choose the maximal parabolic subgroup $P_i$ 
as the stabilizer of $\zeta_i$ in $G$, 
compare the discussion in the introduction 
preceding Theorem \ref{generalschubertinequalitiesintro} 
and the notation used there. 

If the order of the dihedral group is $2m$ 
with $m=3,4$ or $6$
then the complex dimension of each Grassmannian is $m-1$ and there are 
$m$ Bruhat cells, one of each dimension between $0$ and $m-1$. It follows that
the rational cohomology rings are polynomial algebras on a two-dimensional
generator (the hyperplane section class). 

Let $w_i$ be the reflection in $W$ fixing $\zeta_i$ 
and let $\la_{P_i}$ 
be the unique fundamental weight fixed by $w_i$.
The Schubert cycles in $G/P_i$ 
are in one-to-one correspondence 
with the vertices of the spherical Coxeter complex 
in the $W$-orbit of $\zeta_i$,
respectively, 
with the maximally singular weights
in the $W$-orbit of $\la_{P_i}$,
We measure the word length in $W$ with respect to the generating set 
$\{w_1,w_2\}$
and can thus speak of the length of a coset in 
$W^{P_i}:=W/\{e,w_i\}$. 
For $j=0,\dots, m-1$
there is a unique coset of length $j$,
a corresponding weight $\la_j$ in the orbit $W\la_{P_i}$
and a corresponding Schubert cycle $C_{\la_j}$. 
We will abuse notation
and also let $C_{\la_j}$ 
denote the homology class carried by the Schubert cycle 
$C_{\la_j}$. 
Then $C_{\la_j}$ is a generator of the infinite cyclic group 
$H_{2j}(G/P_i;\Z)$. We let
$\ga_{m-j-1}$ be the cohomology class which is the Poincar\'e dual of
$C_{w_j}$. Thus $\ga_j$ is a generator of $H^{2j}(G/P_i;\Z)$. 

The system of stability inequalities divides into 
two subsystems,
one for each Grassmannian. 
The inequalities for each of these  
subsystems are parametrized by a subset of the ordered partitions of $m-1$  
into $n$ nonnegative integers. 
The partition 
$m-1=j_1+\dots+j_n$ 
gives rise to an
inequality in the $G/P_i$-subsystem
if and only if the product $\ga_{j_1} \cdots \ga_{j_n}$ 
is the fundamental class generating
$H^{2m-2}(G/P_i;\Z)$. 
Note that 
the weak triangle inequalities in the sense of section \ref{sec:wtrineq} 
correspond to those decompositions 
$\ga_{m-1}=\ga_j\cdot\ga_k\cdot\ga_l$ of the fundamental class 
where at least one factor has degree zero.
The symmetric group $S_n$ acts naturally on the set of inequalities 
by permuting the $\De$-side lengths.
Thus unordered partitions give rise to $S_n$-orbits of inequalities
and in our examples below we will write down 
one representing inequality from each orbit.

Note that for $m-1=j+k+l$ 
we have $\ga_j \cdot \ga_k = c^i_{jk} \,  \ga_{m-l-1}$
with certain nonnegative integers $c^i_{jk}$ 
which we will refer to as the 
{\em structure constants}
of the ring $H^{\ast}(G/P_i;\Z)$. 
In the case $n=3$
the inequality corresponding to the partition 
$(j,k,l)$ occurs if and only if $c^i_{jk}=1$.
The cohomology rings $H^{\ast}(G/P_i;\Z)$
may be easily calculated using Chevalley's
formula, see Lemma 8.1 of \cite{FultonWoodward} or Theorem 6.1 of 
\cite{TelemanWoodward} as was 
pointed out to us by Chris Woodward. In addition we have taken the cohomology 
rings for $G_2$ from \cite{TelemanWoodward}, page 20.

\subsection{The polyhedron for $A_2$}

We consider the group $G=SL(3,\C)$. 
The Euclidean Weyl chamber is given by 
\begin{equation*}
\Delta_{euc} = \{(x,y,z): x+y+z = 0 \,\land\, x\geq y \geq z \}.
\end{equation*}
The fundamental weights are 
\begin{equation*}
\la_{P_1}(x,y,z)= x \quad \mbox{and} \quad \la_{P_2}(x,y,z)= -z.
\end{equation*}
One Grassmannian is $\C\P^2$ and the other
is  the dual projective space $(\C\P^2)^{\vee}$. 

In $\C\P^2$
the $0$-, $1$- and $2$-dimensional Schubert cycles 
$[pt]$, $[\C\P^1]$ and $[\C\P^2]$
correspond to the maximally singular coweights 
$(2,-1,-1)$, $(-1,2,-1)$ and $(-1,-1,2)$
and hence to the weights 
$x$, $y$ and $z$.
Similarly,
the $0$-, $1$- and $2$-dimensional Schubert cycles 
in $(\C\P^2)^{\vee}$ 
correspond to the maximally singular coweights 
$(1,1,-2)$, $(1,-2,1)$ and $(-2,1,1)$,
respectively,
to the weights 
$-x$, $-y$ and $-z$.
In both cases, 
all the structure constants are $1$ 
and we get one inequality for each unordered partition of $2$. 

As we mentioned before 
the symmetric group $S_n$ acts naturally on the set of inequalities 
by permuting the $\De$-side lengths.
The $S_n$-orbits of inequalities correspond to the ordered partitions of $2$
and we will write down one representing inequality
for each ordered partition. 

\begin{center}
The inequalities associated to $\C\P^2$:
\end{center}
\begin{align*}
x_1 + z_2 + \dots + z_n \leq 0 \\
y_1 + y_2 + z_3 + \dots + z_n \leq 0 
\end{align*}

The two inequalities correspond
to the partitions $2+0+\dots+0$ and $1+1+0+\dots+0$ of $2$,
that is, 
to the one point intersections of Schubert cycles
$[pt]\cdot[\C\P^2]\cdots[\C\P^2]=[pt]$
and 
$[\C\P^1]\cdot[\C\P^1]\cdot[\C\P^2]\cdots[\C\P^2]=[pt]$,
respectively, 
to the decompositions 
$\ga_2\ga_0^{n-1}=\ga_1^2\ga_0^{n-2}$
of the fundamental class $\ga_2$. 
Similarly, we have 

\begin{center}
The inequalities associated to $(\C\P^2)^{\vee}$:
\end{center}
\begin{align*}
z_1 + x_2 + \dots + x_n \geq 0 \\
y_1 + y_2 + x_3 + \dots + x_n \geq 0 
\end{align*}

In the case $n=3$
all inequalities are weak triangle inequalities 
and moreover
the system of inequalities is known to be irredundant, 
see \cite{KnutsonTaoWoodward}.

The following $8$ vectors are a set of generators of the polyhedral cone 
${\cal P}_3(A_2)$ in the $6$-dimensional space ${\goth a}^3$. 
\begin{center}
\setlength{\tabcolsep}{ .75cm}
\renewcommand{\arraystretch}{1.5}

\begin{tabular}{l c l c l c}
(2,-1,-1) (2,-1,-1) (2,-1,-1) &   (1,1,-2)  (2,-1,-1) (0,0,0)  \\
(0,0,0)   (1,1,-2) (2,-1,-1)  &   (2,-1,-1) (0,0,0)   (1,1,-2)  \\ 
(0,0,0)   (2,-1,-1) (1,1,-2)  &   (2,-1,-1) (1,1,-2)  (0,0,0)   \\ 
(1,1,-2)   (0,0,0)   (2,-1,-1)  &   (1,1,-2)    (1,1,-2)  (1,1,-2)
\end{tabular}
\end{center}

\subsection{The polyhedron for $B_2=C_2$.}

We consider the group $G=SO(5)$.
(Note that $SO(5)$ is isomorphic to $PSp(4)$.)

After a suitable change of coordinates 
we may write the invariant quadratic form 
on $\C^5$ as 
$q(u)=u_1u_5+u_2u_4+u_3^2$. 
The diagonal matrices with real eigenvalues 
$x,y,0,-x,-y$ 
then form a Cartan subalgebra ${\goth a}$ in ${\goth g}$ 
and the Euclidean Weyl chamber is given by
\begin{equation*}
\Delta_{euc} = \{(x,y):  x\geq y \geq 0 \}.
\end{equation*}
The fundamental weights are 
\begin{equation*}
\la_{P_1}(x,y)= x \quad \mbox{and} \quad \la_{P_2}(x,y)= x+y.
\end{equation*}

The Grassmannian $G/P_1$ is the space of isotropic lines 
in $\C^5$. 
This is the smooth quadric three-fold $\mathcal{Q}_3$ in $\C\P^4$
given by the equation $q(u)=0$. 
The other Grassmannian $G/P_2$ is the space of 
totally-isotropic two-planes in $\C^5$. 

For the Grassmannian $G/P_1$
the Schubert cycles of dimension $0,1,2$ and $3$
correspond to the maximally singular coweights 
$(1,0),(0,1),(0,-1)$ and $(-1,0)$, respectively,
to the weights
$\la=x,y,-y$ and $-x$:

\begin{center}
\begin{tabular} {|l|c|l|}
\hline
$\lambda$ & dim $C_{\lambda}$ & PD $C_{\lambda}$ \\ \hline
   $x$     & $0$ & $\ga_3$  \\ \hline
   $y$     & $1$ & $\ga_2$  \\ \hline
   $-y$    & $2$ & $\ga_1$   \\ \hline
   $-x$    & $3$ & $1$     \\ \hline
\end{tabular}
\end{center}

The Schubert cycles and their homological intersections 
can still be determined by hand: 
The $2$-dimensional Schubert cycle is a hyperplane section
and the $1$-dimensional Schubert cycle is an embedded projective line. 
In particular,
the self intersection of the $2$-cycle is {\em twice}
the $1$-cycle. 
The cohomology ring is given by the following table:

\begin{center}
\begin{tabular} {|c|l|l|l|l|}
\hline
$H^*({\cal Q}_3)$ & 1 & $\ga_1$ & $\ga_2$ & $\ga_3$ \\ \hline
   1         & 1 & $\ga_1$ & $\ga_2$ & $\ga_3$ \\ \hline
   $\ga_1$     & $\ga_1$ & $2\ga_2$ & $\ga_3$ & $0$ \\ \hline
   $\ga_2$     & $\ga_2$ & $\ga_3$  & $0$ & $0$ \\ \hline
   $\ga_3$     & $\ga_3$ & $0$    & $0$ & $0$ \\ \hline
\end{tabular}
\end{center}
 
We see that the structure constant $c^1_{1,1}$ equals $2$,
so any partition of 3 which involves the pair $(1,1)$ 
does not give rise to an inequality.
Indeed, the only decompositions of the fundamental class 
into products of Schubert classes are 
$\ga_3=\ga_2\cdot\ga_1$.

\begin{center}
The inequalities associated to $G/P_1={\cal Q}_3$:
\end{center}

\begin{align*}
x_1 \leq x_2 + \dots + x_n \\
y_1 - y_2 \leq x_3 + \dots + x_n 
\end{align*}

For the Grassmannian $G/P_2$
the Schubert cycles of dimension $0,1,2$ and $3$
correspond to the maximally singular coweights 
$(1,1),(1,-1),(-1,1)$ and $(-1,-1)$, respectively,
to the weights
$\la=x+y,x-y,-x+y$ and $-x-y$:

\begin{center}
\begin{tabular} {|l|c|l|}
\hline
$\lambda$ & dim $C_{\lambda}$ & PD $C_{\lambda}$ \\ \hline
   $x+y$     & $0$ & $\ga_3$  \\ \hline
   $x-y$    & $1$ & $\ga_2$  \\ \hline
   $-x+y$    & $2$ & $\ga_1$   \\ \hline
   $-x-y$   & $3$ & $1$     \\ \hline
\end{tabular}
\end{center}

The cohomology ring of $G/P_2$ is easily determined
due to the exceptional isomorphism 
$SO(5,\C)\cong PSp(4,C)$.
(Recall that $Sp(4,\C)$, 
the automorphism group of a complex symplectic 2-form $\om$ on $\C^4$,
acts on the $6$-dimensional vector space 
$\La^2(\C^4)^{\ast}$ 
of alternating bilinear forms on $\C^4$. 
The induced action 
$Sp(4,\C)\acts\La^2(\C^4)^{\ast}$ 
preserves a natural non-degenerate quadratic form 
and the line generated by $\om$.)
The Grassmannians associated to the groups 
$SO(5,\C)$ and $Sp(4,C)$ 
are the same and one of the Grassmannians for $Sp(4,C)$ is $\C\P^3$.
Since $G/P_1\cong{\cal Q}_3\not\cong\C\P^3$
we conclude that 
$G/P_2\cong\C\P^3$.
Thus the structure constants for the cohomology ring of $G/P_2$ 
are given by the following table:
 
\begin{center}
\begin{tabular} {|c|l|l|l|l|}
\hline
$H^*(G/P_2)$ & 1 & $\ga_1$ & $\ga_2$ & $\ga_3$ \\ \hline
   1         & 1 & $\ga_1$ & $\ga_2$ & $\ga_3$ \\ \hline
   $\ga_1$     & $\ga_1$ & $\ga_2$ & $\ga_3$ & $0$ \\ \hline
   $\ga_2$     & $\ga_2$ & $\ga_3$  & $0$ & $0$ \\ \hline
   $\ga_3$     & $\ga_3$ & $0$    & $0$ & $0$ \\ \hline
\end{tabular}
\end{center}

The possible decompositions of the fundamental class 
into products of Schubert classes are 
$\ga_3=\ga_2\cdot\ga_1=\ga_1^3$.

\begin{center}
The inequalities associated to $G/P_2$:
\end{center}

\begin{align*}
x_1 + y_1 \leq x_2 + y_2 + \dots + x_n + y_n \\
x_1 - y_1 -x_2 + y_2  \leq x_3 + y_3 + \dots + x_n + y_n  \\
-x_1 + y_1 - x_2 + y_2 -x_3 + y_3 \leq x_4 + y_4 + \dots + x_n + y_n 
\end{align*}

Observe that this last inequality is redundant
because it follows from $y\leq x$
which is one of the defining inequalities of $\De_{euc}$. 

In the case $n=3$,
according to Porta the system obtained when the last 
($S_3$-invariant) inequality is
removed is irredundant.
We observe that the redundant triangle inequalities 
are exactly the non-weak ones.

The following $12$ vectors are a set of generators of the polyhedral cone 
${\cal P}_3(B_2) = {\cal P}_3(C_2)$ 
in the $6$-dimensional space ${\goth a}^3$. 
\begin{center}
\setlength{\tabcolsep}{ .75cm}
\renewcommand{\arraystretch}{1.5}

\begin{tabular}{l c l c l c}
(1,1) (1,1) (2,0)   & (1,0) (0,0) (1,0)   & (1,1) (1,1) (0,0) \\
(1,1) (2,0) (1,1)   & (1,0) (1,0) (0,0)   & (1,0) (1,0) (1,1) \\
(2,0) (1,1) (1,1)   & (0,0) (1,1) (1,1)   & (1,0) (1,1) (1,0) \\
(0,0) (1,0) (1,0)   & (1,1) (0,0) (1,1)   & (1,1) (1,0) (1,0) 
\end{tabular}
\end{center}

\subsection{The polyhedron for $G_2$} 

In this case both Grassmannians have dimension $5$. 

We use non-rectangular linear coordinates on ${\goth a}$ 
such that the Euclidean Weyl chamber is given by 
\begin{equation*}
\De_{euc} = \{(x,y): x\geq 0 , y\geq 0\} 
\end{equation*}
and such that the standard basis vectors are the fundamental coweights.
We require moreover that 
$\|(1,0)\|=\sqrt{3}\cdot\|(0,1)\|$.
With respect to these coordinates 
the natural Euclidean metric on ${\goth a}$ takes, up to scale, the form 
$3dx\otimes dx+\frac{3}{2}dx\otimes dy+\frac{3}{2}dy\otimes dx+dy\otimes dy$.

Let $G/P_1$ be the Grassmannian 
corresponding to the fundamental coweight $(1,0)$.
The Schubert cycles of dimension $0,\dots, 5$
correspond to the Weyl group orbit of maximally singular coweights 
$(1,0)$, $(-1,3)$, $(2,-3)$, $(-2,3)$, $(1,-3)$, $(-1,0)$, 
respectively,
via the scalar product on ${\goth a}$ up to a scale factor
to the orbit of weights
$2x+y$, $x+y$, $x$, $-x$, $-x-y$, $-2x-y$.

\begin{center}
\begin{tabular} {|l|c|l|}
\hline
$\lambda$ & dim $C_{\lambda}$ & PD $C_{\lambda}$ \\ \hline
   $2x+y$    & $0$ & $\ga_5$   \\ \hline
   $x+y$    & $1$ & $\ga_4$   \\ \hline
   $x$    & $2$ & $\ga_3$   \\ \hline
   $-x$   & $3$ & $\ga_2$   \\ \hline
   $-x-y$  & $4$ & $\ga_1$   \\ \hline
   $-2x-y$  & $5$ & $1$     \\ \hline
\end{tabular}
\end{center}

The structure constants for the cohomology ring of $G/P_1$ 
are given by the following table \cite{TelemanWoodward}: 

\begin{center}
\begin{tabular} {|c|l|l|l|l|l|l|}
\hline
$H^*(G/P_1)$ & 1     & $\ga_1$  & $\ga_2$  & $\ga_3$ & $\ga_4$ & $\ga_5$ \\ \hline
   1         & 1     & $\ga_1$  & $\ga_2$  & $\ga_3$ & $\ga_4$ & $\ga_5$ \\ \hline
   $\ga_1$     & $\ga_1$ & $\ga_2$  & $2\ga_3$ & $\ga_4$ & $\ga_5$ & $0$ \\ \hline
   $\ga_2$     & $\ga_2$ & $2\ga_3$ & $2\ga_4$ & $\ga_5$ & $0$   & $0$ \\ \hline
   $\ga_3$     & $\ga_3$ & $\ga_4$  & $\ga_5$  & $0$   & $0$   & $0$ \\ \hline
   $\ga_4$     & $\ga_4$ & $\ga_5$  & $0$    & $0$   & $0$   & $0$ \\ \hline
   $\ga_5$     & $\ga_5$ & $0$    & $0$    & $0$   & $0$   & $0$ \\ \hline
\end{tabular}
\end{center}

The possible decompositions of the fundamental class 
into products of Schubert classes are 
$\ga_5=\ga_4\cdot\ga_1=\ga_3\cdot\ga_2=\ga_3\cdot\ga_1^2$.

\begin{center}
The inequalities associated to $G/P_1$:
\end{center}

\begin{align*}
2x_1 + y_1 \leq 2x_2 + y_2 + \dots + 2x_n + y_n \\
x_1 + y_1 - x_2 - y_2  \leq 2x_3 + y_3 + \dots + 2x_n + y_n \\
x_1 - x_2 \leq 2x_3 + y_3 + \dots + 2x_n + y_n \\ 
x_1 - x_2 - y_2  - x_3 - y_3  \leq 2x_4 + y_4 + \dots + 2x_n + y_n
\end{align*}

Note that the fourth group of inequalities is redundant.
Indeed, we may obtain it from the first inequality 
by adding the inequality 
$0\leq y_1+y_2+y_3+2x_4 + y_4 + \dots + 2x_n + y_n$
which is implied by the inequalities defining $\De_{euc}$.

We now describe the subsystem corresponding to $G/P_2$,
the Grassmannian 
corresponding to the fundamental coweight $(0,1)$.
The Schubert cycles of dimension $0,\dots,5$
correspond to the Weyl group orbit of maximally singular coweights 
$(0,1)$, $(1,-1)$, $(-1,2)$, $(1,-2)$, $(-1,1)$, $(0,-1)$,
respectively, 
to the orbit of weights 
$3x+2y$, $3x+y$, $y$, $-y$, $-3x-y$, $-3x-2y$.

\begin{center}
\begin{tabular} {|l|c|l|}
\hline
$\lambda$ & dim $C_{\lambda}$ & PD $C_{\lambda}$ \\ \hline
   $3x+2y$    & $0$ & $\ga_5$   \\ \hline
   $3x+y$    & $1$ & $\ga_4$   \\ \hline
   $y$    & $2$ & $\ga_3$   \\ \hline
   $-y$   & $3$ & $\ga_2$   \\ \hline
   $-3x-y$  & $4$ & $\ga_1$   \\ \hline
   $-3x-2y$  & $5$ & $1$     \\ \hline
\end{tabular}
\end{center}

The structure constants for the cohomology ring of $G/P_2$ 
are given by the following table \cite{TelemanWoodward}:

\begin{center}
\begin{tabular} {|c|l|l|l|l|l|l|}
\hline
$H^*(G/P_2)$ & 1     & $\ga_1$  & $\ga_2$  & $\ga_3$ & $\ga_4$ & $\ga_5$ \\ \hline
   1         & 1     & $\ga_1$  & $\ga_2$  & $\ga_3$ & $\ga_4$ & $\ga_5$ \\ \hline
   $\ga_1$     & $\ga_1$ & $3\ga_2$ & $2\ga_3$ & $3\ga_4$& $\ga_5$ & $0$ \\ \hline
   $\ga_2$     & $\ga_2$ & $2\ga_3$ & $2\ga_4$ & $\ga_5$ & $0$   & $0$ \\ \hline
   $\ga_3$     & $\ga_3$ & $3\ga_4$ & $\ga_5$  & $0$   & $0$   & $0$ \\ \hline
   $\ga_4$     & $\ga_4$ & $\ga_5$  & $0$    & $0$   & $0$   & $0$ \\ \hline
   $\ga_5$     & $\ga_5$ & $0$    & $0$    & $0$   & $0$   & $0$ \\ \hline
\end{tabular}
\end{center}

The possible decompositions of the fundamental class 
into products of Schubert classes are 
$\ga_5=\ga_4\cdot\ga_1=\ga_3\cdot\ga_2$. 

\begin{center}
The inequalities associated to $G/P_2$:
\end{center}

\begin{align*}
3x_1 + 2y_1 \leq 3x_2 + 2y_2 + \dots + 3x_n + 2y_n \\ 
3x_1 + y_1 - 3x_2 - y_2 \leq 3x_3 + 2y_3 + \dots + 3x_n + 2y_n \\ 
y_1 - y_2 \leq 3x_3 + 2y_3 + \dots + 3x_n + 2y_n 
\end{align*}

Our system of inequalities 
becomes that of \cite[p.\ 458]{BerensteinSjamaar} once one replaces
$y_j$ in our inequalities with $3 y_j$ and takes into account that
they have extra inequalities corresponding to intersections of Schubert
classes that are nonzero multiples of the point class 
that are not equal to $1$.

In the case $n=3$
we find using the program Porta 
that the system obtained by
removing the $S_3$-orbit of redundant inequalities mentioned previously 
is an irredundant system.
We observe that also in this case 
the redundant triangle inequalities are exactly the non-weak ones. 

The following $24$ vectors are a set of generators of the polyhedral cone 
${\cal P}_3(G_2)$ in the $6$-dimensional space ${\goth a}^3$.
\begin{center}
\setlength{\tabcolsep}{ .75cm}
\renewcommand{\arraystretch}{1.5}

\begin{tabular}{l c l c l c}
(0,3) (1,0) (2,0)  &  (1,0) (0,1) (0,2)  &  (0,1) (0,0) (0,1) \\
(0,3) (2,0) (1,0)  &  (1,0) (0,2) (0,1)  &  (0,1) (0,1) (0,0) \\
(1,0) (0,3) (2,0)  &  (0,2) (0,1) (1,0)  &  (0,0) (1,0) (1,0) \\
(1,0) (2,0) (0,3)  &  (0,2) (1,0) (0,1)  &  (1,0) (0,0) (1,0) \\
(2,0) (0,3) (1,0)  &  (0,3) (1,0) (1,0)  &  (1,0) (1,0) (0,0) \\
(2,0) (1,0) (0,3)  &  (1,0) (0,3) (1,0)  &  (0,1) (0,1) (1,0) \\
(0,1) (0,2) (1,0)  &  (1,0) (1,0) (0,3)  &  (0,1) (1,0) (0,1)  \\
(0,1) (1,0) (0,2)  &  (0,0) (0,1) (0,1)  &  (1,0) (0,1) (0,1)  
\end{tabular}
\end{center}

\medskip\no
Michael Kapovich, 
Department of Mathematics, 
University of California, 
Davis, CA 95616, USA,
kapovich@math.ucdavis.edu

\no
Bernhard Leeb, 
Mathematisches Institut,
Universit\"at M\"unchen, 
Theresienstrasse 39,
D-80333 M\"unchen, Germany, 
b.l@lmu.de

\no
John Millson, 
Department of Mathematics, 
University of Maryland, 
College Park, MD 20742, USA, 
jjm@math.umd.edu

\end{document}